\documentclass{gtart_h}

\def\ifplaintex{\expandafter\ifx\csname documentclass\endcsname\relax}

\ifplaintex 
\else
\fi

\def\gt{{\mathsurround=0pt\it $\cal G\mskip-2mu$eometry \&\ 
$\cal T\!\!$opology}}        %  journal title in recommended style

\def\gtm{{\mathsurround=0pt\it $\cal G\mskip-2mu$eometry \&\ 
$\cal T\!\!$opology $\cal M\mskip-1mu$onographs}}    %  for monographs

\def\gtp{{\mathsurround=0pt\it $\cal G\mskip-2mu$eometry \&\ 
$\cal T\!\!$opology $\cal P\!$ublications}}  % GT publications

\def\recd{{\small Received:\qua\receiveddate\ifx\reviseddate\relax
\else\qquad Revised:\qua\reviseddate\fi\par}} 

\def\volumenumber#1{\def\thevolumenumber{#1}}
\def\volumeyear#1{\def\thevolumeyear{#1}}
\def\volumename#1{\def\thevolumename{#1}}
\def\papernumber#1{\def\thepapernumber{#1}}
\def\pagenumbers#1#2{\def\startpage{#1}\def\finishpage{#2}}
\def\published#1{\def\publishdate{#1}}
\def\received#1{\def\receiveddate{#1}}
\def\revised#1{\def\reviseddate{#1}}
\def\accepted#1{\def\accepteddate{#1}}

\def\asciiaddress#1{\def\theasciiaddress{#1}}
\def\asciiemail#1{\def\theasciiemail{#1}}
\def\asciiurl#1{\def\theasciiurl{#1}}

\long\def\asciiabstract#1{\long\def\theasciiabstract{#1}}

\let\\\par
\let\thevolumenumber\relax\let\thepapernumber\relax
\let\thevolumeyear\relax\let\startpage\relax
\let\finishpage\relax\let\publishdate\relax\let\receiveddate\relax
\let\reviseddate\relax\let\accepteddate\relax\let\theasciititle\relax
\let\theasciiauthors\relax\let\theasciiaddress\relax
\let\theasciiabstract\relax

\let\theerratum\relax\let\theasciiemail\relax
\let\theshortauthors\relax\let\theshorttitle\relax\let\theasciiurl\relax

\def\startpage{1}\def\finishpage{15}\def\thepapernumber{77}

\volumenumber{2}
\volumename{Proceedings of the Kirbyfest}
\volumeyear{1999}

\long\def\maketitlep{   % start of definition of \maketitlep

\count0=\startpage

\gtm\nl        %   GT mongraphs (top left) 
{\small Volume \thevolumenumber: \thevolumename\nl 
\ifx\theerratum\relax\else Erratum \erratumnumber\nl\fi
Pages \startpage--\finishpage\nl}

\vglue 0.1truein   % top margin

{\parskip=0pt\leftskip 0pt plus 1fil\def\\{\par\smallskip}{\ifplaintex\large
\else\Large\fi\bf\thetitle}\par\medskip}   
\vglue 0.05truein 

{\parskip=0pt\leftskip 0pt plus 1fil\def\\{\par}{\sc\theauthors}
\par\medskip}%
 
\vglue 0.03truein 

{\small\leftskip 25pt\rightskip 25pt{\bf Abstract}\stdspace\theabstract

{\bf AMS Classification}\stdspace\theprimaryclass
\ifx\thesecondaryclass\relax\else; \thesecondaryclass\fi\par
{\bf Keywords}\stdspace \thekeywords\par}\vglue 7pt

}   % end of definition of \maketitlep

\font\phead=cmsl9 scaled 950
\font=cmsl9 scaled 1050
\font\pnum=cmbx10 scaled 913
\font\lnum=cmbx10 
\font\pfoot=cmsl9 scaled 950
\font=cmsl9 scaled 1050
\ifplaintex
\headline{\vbox to 0pt{\vskip -4.5mm\line{\small\phead\ifnum
\count0=\startpage ISSN 1464-8997 (on line)
1464-8989 (printed) \hfill {\pnum\folio}\else\ifodd\count0\def\\{ }% 
\ifx\theshorttitle\relax\thetitle\else\theshorttitle\fi\hfill{\pnum\folio}
\else\def\\{ and }{\pnum\folio}\hfill\ifx\theshortauthors\relax\theauthors
\else\theshortauthors\fi\fi\fi}\vss}}
\footline{\vbox to 0pt{\vglue 0mm\line{\small\pfoot\ifnum\count0=\startpage
Published \publishdate:\qua\copyright\ \gtp\hfill\else
\gtm, Volume \thevolumenumber\ (\thevolumeyear)\hfill\fi}\vss
}}
\else
\makeatletter
\def\@oddhead{{\small\lhead\ifnum\count0=\startpage ISSN 1464-8997 (on line)
1464-8989 (printed) \hfill {\lnum\number\count0}\else\ifodd\count0
\def\\{ }\ifx\theshorttitle\relax \thetitle \else\theshorttitle\fi\hfill
{\lnum\number\count0}\else\def\\{ and }{\lnum\number\count0}
\hfill\ifx\theshortauthors\relax 
\theauthors\else\theshortauthors\fi\fi\fi}}\def\@evenhead{@oddhead}
\def\@oddfoot{\small\lfoot\ifnum\count0=\startpage Published \publishdate:\qua\copyright\ \gtp\hfill\else
\gtm, Volume \thevolumenumber\ (\thevolumeyear)\hfill\fi}
\def\@evenfoot{@oddfoot}
\makeatother
\fi

\let\maketitlepage\maketitlep
\let\makeshorttitle\maketitlepage
\let\maketitle\maketitlepage

\newwrite\gtoutfile
\long\gdef\makeheadfile{  %%% start of definition of \makeheadfile
{\def\\{, }\def\s{ }
\immediate\openout\gtoutfile head.xxx
\immediate\write\gtoutfile{Proxy-for: \ifx\theasciiauthors\relax
\theauthors\else\theasciiauthors\fi\s<\ifx\theasciiemail\relax\theemail\else\theasciiemail\fi>}
\immediate\write\gtoutfile{\noexpand\\}
\immediate\write\gtoutfile{Authors: \ifx\theasciiauthors\relax
\theauthors\else\theasciiauthors\fi}
{\def\\{ }\immediate\write\gtoutfile{Title: \ifx\theasciititle\relax
\thetitle\else\theasciititle\fi}}
\immediate\write\gtoutfile{Subj-class: GT or SG, GR etc}
\immediate\write\gtoutfile{MSC-class: \theprimaryclass\ifx\thesecondaryclass\relax\else, \thesecondaryclass\fi}
\immediate\write\gtoutfile{Journal-ref: Geom. Topol. Monogr. \thevolumenumber\s
(\thevolumeyear) \startpage-\finishpage}
\immediate\write\gtoutfile{Comments: Published by Geometry and Topology Monographs at}
\immediate\write\gtoutfile{\s\s\s  http://www.maths.warwick.ac.uk/gt/GTMon\thevolumenumber/paper\thepapernumber.abs.html}
\immediate\write\gtoutfile{\noexpand\\}
\immediate\write\gtoutfile{}
\ifx\theasciiabstract\relax
\immediate\write\gtoutfile{\theabstract}\else
\immediate\write\gtoutfile{\theasciiabstract}\fi
\immediate\write\gtoutfile{}
\immediate\write\gtoutfile{\noexpand\\}
\immediate\write\gtoutfile{}
\immediate\closeout\gtoutfile}}  %%% end of definition of \makeheadfile

\def\maketitlepage{\maketitlep\makeheadfile}
\let\makeshorttitle\maketitlepage
\let\maketitle\maketitlepage

\volumenumber{7}
\volumename{Proceedings of the Casson Fest}
\volumeyear{2004}
\papernumber{17}
\pagenumbers{509}{547}
\received{6 October 2003}
\revised{21 April 2005}
\accepted{21 April 2005}
\published{21 May 2005}
\usepackage{amsmath,amssymb,graphicx}

\def\leq{\leqslant}
\def\geq{\geqslant}

\newtheorem{theorem}{Theorem}
\newtheorem{proposition}[theorem]{Proposition}
\newtheorem{lemma}[theorem]{Lemma}

\newtheorem{corollary}[theorem]{Corollary}

\theoremstyle{remark}
\newtheorem{remark}[theorem]{Remark}

\begin{document}

\title [The metric space of geodesic laminations]{The metric
space of geodesic laminations\\on a surface II: small surfaces}

\authors{Francis Bonahon\\Xiaodong Zhu}

\address{Department of Mathematics,  University of Southern
California\\Los Angeles, CA~90089--2532, USA}
\gtemail{\mailto{fbonahon@math.usc.edu}, \mailto{xzhu@juniper.net}}
\asciiemail{fbonahon@math.usc.edu, xzhu@juniper.net}

\secondaddress{Juniper Networks, 1194 North Mathilda Avenue\\Sunnyvale,
CA 94089--1206, USA}
\urladdr{http://www-rcf.usc.edu/~fbonahon/}
\asciiurl{http://www-rcf.usc.edu/ fbonahon/}

\asciiaddress{Department of Mathematics,  University of Southern
California\\Los Angeles, CA 90089-2532, USA\\and\\Juniper Networks, 
1194 North Mathilda Avenue\\Sunnyvale, CA 94089-1206, USA}

\begin{abstract} 
We continue our investigation of the space of geodesic laminations on
a surface, endowed with the Hausdorff topology. We determine the
topology of this space for the once-punctured torus and the
4--times-punctured sphere. For these two surfaces, we also compute the
Hausdorff dimension of the space of geodesic laminations, when it is
endowed with the natural metric which, for small distances, is $-1$
over the logarithm of the Hausdorff metric. The key ingredient is an
estimate of the Hausdorff metric between two simple closed geodesics
in terms of their respective slopes.
\end{abstract}

\asciiabstract{%
We continue our investigation of the space of geodesic laminations on
a surface, endowed with the Hausdorff topology. We determine the
topology of this space for the once-punctured torus and the
4-times-punctured sphere. For these two surfaces, we also compute the
Hausdorff dimension of the space of geodesic laminations, when it is
endowed with the natural metric which, for small distances, is -1
over the logarithm of the Hausdorff metric. The key ingredient is an
estimate of the Hausdorff metric between two simple closed geodesics
in terms of their respective slopes.}


\maketitle

This article is a continuation of
 the study of the
Hausdorff metric $d_{\mathrm
H}$ on
the space $\mathcal L(S)$ of all
geodesic laminations on a
surface $S$, which we began in the
article \cite{ZhuBon1}. The
impetus for these two papers
originated in  the
monograph
\cite{CasBle} by Andrew Casson
and Steve Bleiler,  which was the
first to systematically exploit
the Hausdorff topology on the
space of geodesic laminations.

In this paper, we
restrict  attention to the
case where the surface $S$ is the
once-punctured torus or the
4--times-punctured sphere. To
some extent, these are the first
non-trivial examples, since
$\mathcal L(S)$ is defined
only when the Euler
characteristic of $S$ is
negative, is finite when
$S$ is the 3--times-punctured
sphere or the twice-punctured
projective plane, and is countable
infinite when $S$ is the
once-punctured Klein bottle (see
for instance Section~\ref
{sect:VerySmall}).

We
will also restrict attention to
the open and closed subset
$\mathcal L_0(S)$ of
$\mathcal L(S)$ consisting of
those geodesic laminations which
are disjoint from the boundary.
This second restriction is
only an expository choice. The
results and techniques of the
paper can be relatively easily 
extended to the full space
$\mathcal L(S)$, but at the
expense of many more cases to
consider; the corresponding
strengthening of the
results did not seem to be worth
the increase in size of the
article. 

The first two results deal with
the topology of $\mathcal L_0(S)$
for these two surfaces.

\begin{theorem}
\label{thm:TopPunctTor}
When $S$ is the once-punctured
torus, the space $\mathcal
L_0(S)$ naturally splits as the
disjoint union of two compact
subsets, the closure $\mathcal
L_0^{\mathrm{cr}}(S)$ of the
set of simple closed curves and
its complement
$\mathcal L_0(S)-\mathcal
L_0^{\mathrm{cr}}(S)$. The first
subspace $\mathcal
L_0^{\mathrm{cr}}(S)$ is
homeomorphic to a subspace $K\cup
L_1$ of the circle $\mathbb S^1$,
where $K$ is the
standard Cantor set and where $L_1$
is a countable set consisting of one
isolated point in each component
of
$\mathbb S^1-K$. The complement 
$\mathcal L_0(S)-\mathcal
L_0^{\mathrm{cr}}(S)$ is
homeomorphic to a subspace $K\cup
L_3$ of
$\mathbb S^1$, union of the Cantor
set
$K\subset \mathbb S^1$ and of a
countable set $L_3$ consisting of
exactly  $3$ isolated points in
each component of
$\mathbb S^1-K$.  
\end{theorem}

\begin{theorem} 
\label{thm:TopPunctSpher}
When $S$ is the
$4$--times-punctured sphere, 
the space $\mathcal
L_0(S)$ is
homeomorphic to a subspace $K\cup
L_7$ of
$\mathbb S^1$, union of the Cantor
set
$K$ and of a
countable set $L_7$ consisting of
exactly  $7$ isolated points in
each component of
$\mathbb S^1-K$. In this case,
the closure $\mathcal
L^{\mathrm{cr}}_0(S)$ of the set
 of simple closed curves
is the union $K\cup L_1$ of $K$ and
of a discrete set
$L_1\subset L_7$ consisting of
exactly one point in each component
of
$\mathbb S^1-K$; in particular, its
complement 
$\mathcal L_0(S)-\mathcal
L_0^{\mathrm{cr}}(S)$
is countable infinite. 
\end{theorem}

The above subspaces $K\cup L_1$, $K\cup L_3$ and 
$K\cup L_7$ are all homeomorphic. 
However, it is convenient to keep a
distinction between these spaces,
because the proofs of
Theorems~\ref{thm:TopPunctTor} and
\ref{thm:TopPunctSpher} make the
corresponding embeddings of $\mathcal
L_0(S)$ and $\mathcal
L^{\mathrm{cr}}_0(S)$ in $\mathbb
S^1$ relatively natural. In
particular, these establish a
one-to-one correspondence between
the components of $\mathbb S^1-K$
and the simple closed curves of
$S$. These embeddings are also well
behaved with respect to the action
of the homeomorphism group of $S$
on $\mathcal
L_0(S)$.

We now consider metric properties
of the Hausdorff metric
$d_{\mathrm H}$ on $\mathcal
L_0(S)$. In
\cite{ZhuBon1}, we showed that
the metric space
$\left(
\mathcal L(S), d_{\mathrm H}
\right)$ has Hausdorff dimension
0. In particular, it is totally
disconnected, which is consistent
with
Theorems~\ref{thm:TopPunctTor}
and~\ref{thm:TopPunctSpher}.
However, we also observed that, to
some extent, the Hausdorff metric
$d_{\mathrm H} $ of $\mathcal
L(S)$ is not very canonical
because it is only defined up to
H\"older equivalence. This lead
us to consider  on $\mathcal
L(S)$ another metric
$d_{\log}$ which, for small
distances, is just equal to
$-1/\log{d_{\mathrm H}}$. This
new metric  $d_{\log}$  has better
invariance properties because it
is well-defined up to Lipschitz
equivalence; in particular, its
Hausdorff dimension is
well-defined. We refer to
\cite{ZhuBon1} and
Section~\ref{sect:TrainTracks} for
precise definitions. 

\begin{theorem} 
\label{thm:HausDim}
When $S$ is the
once-punctured torus or the
$4$--times-punctured sphere, the
Hausdorff dimension of the metric
space $\left (\mathcal L_0(S),
d_{\log} \right)$ is equal to
$2$. Its $2$--dimensional
Hausdorff measure is equal to
$0$. 
\end{theorem}

Theorem~\ref{thm:HausDim} was
used in
\cite{ZhuBon1} to show that, for
a general surface $S$ of negative
Euler characteristic which is not
the 3--times-punctured sphere,
the twice-punctured projective
plane or the once-punctured Klein
bottle, the Hausdorff dimension of
$\left (\mathcal L_0(S), d_{\log}
\right)$ is positive and
finite. 

These results should be contrasted
with the more familiar
Thurston completion of the set of
simple closed curves on $S$, by the
space $\mathcal{PML}(S)$ of
projective measured laminations
\cite{FatLauPoe, PenHar}. For the
once--punctured torus and the 4--times-punctured sphere,
$\mathcal{PML}(S)$ is homeomorphic
to the circle and has Hausdorff
dimension 1.

What
is special about the once-punctured torus and the 4--times-punctured sphere is that there is
a relatively simple classification
of their simple closed curves, or
more generally of their recurrent
geodesic laminations, in terms of
their slope. The key technical
result of this article is an
estimate, proved in Section~\ref
{sect:EstimatesPunctTorus},
which relates the Hausdorff
distance of two simple closed
curves to their slopes. 

\begin{proposition}
Let $\lambda$ and $\lambda'$ be
two simple closed geodesics on
the once-punctured torus or on the
$4$--times-punctured sphere, with
respective slopes $\frac pq <\frac
{p'}{q'} \in \mathbb Q
\cup \{\infty\}$. Their
Hausdorff distance $d_{\mathrm
H}(\lambda, \lambda')$ is such
that 
\begin{equation*}
\mathrm e ^{-c_1 / d\bigl ( \frac
pq,
\frac {p'}{q'} \bigr)}\leq
d_{\mathrm H}(\lambda, \lambda')
\leq
\mathrm e ^{-c_2 / d\bigl ( \frac
pq,
\frac {p'}{q'} \bigr)}
\end{equation*}
where the constants $c_1$,
$c_2>0$ depend only on the metric
on the surface, and where 
\begin{equation*} 
\textstyle 
d\bigl(
\frac pq,
\frac{p'}{q'}\bigr) =\max
\Bigl\{  \frac 1{\left\vert
p'' \right\vert +\left\vert
q''\right\vert};
\frac pq
\leq
\frac{p''}{q''} \leq
\frac{p'}{q'} \Bigr \}.
\end{equation*}
\end{proposition}

The other key ingredient is an
analysis of the above metric $d$
on
$\mathbb Q \cup \{\infty\}$,
which is provided in the
Appendix. 

A large number of the results of
this paper were part of the
dissertation \cite{Zhu}. 

\textbf{Acknowledgements}\qua
{The
two authors were greatly
influenced by Andrew Casson, the
first one directly, the second
one indirectly.  It
is a pleasure to acknowledge our
debt to his work, and to his
personal influence over several
generations of topologists. 

The authors are also very grateful to the referee for a critical
reading of the first version of this article. This work was
partially supported by grants DMS-9504282, DMS-9803445 and
DMS-0103511 from the National Science Foundation.}

\section{Train tracks}
\label{sect:TrainTracks}
We will not repeat the basic
definitions on geodesic
laminations, referring instead to
the standard literature
\cite{CasBle,
CanEpsGre, PenHar, Bon10}, or
to
\cite{ZhuBon1}. However, it is
probably worth reminding the
reader of our definition of the
\emph{Hausdorff distance}
$d_{\mathrm H}(\lambda,
\lambda')$ between two geodesic
laminations $\lambda$,
$\lambda'$ on the surface $S$,
namely 
\begin{equation*}
d_{\mathrm H}(\lambda, \lambda')=
\min
\left\{\varepsilon; \ 
\begin{aligned}
\forall x\in\lambda, \exists
x'\in \lambda' ,d\left((x,
T_x\lambda), (x',
T_{x'}\lambda')\right)
<\varepsilon\\
\forall x'\in\lambda', \exists
x\in \lambda, d\left((x,
T_x\lambda), (x',
T_{x'}\lambda')\right)
<\varepsilon
\end{aligned}
\right\}
\end{equation*}
where the distance $d$ is
measured in the projective
tangent bundle $PT(S)$ consisting
of all pairs $(x, l)$ with $x
\in S$ and with $l$  a line
through the origin in the tangent
space
$T_xS$, and where $T_x\lambda$
denotes the tangent line at $x$ of
the leaf of $\lambda$ passing
through $x$. In particular, $d_{\mathrm H}(\lambda,
\lambda')$ is not the Hausdorff
distance between $\lambda$ and
$\lambda'$ considered as
closed subsets of $S$, but the
Hausdorff distance between their
canonical lifts to $PT(S)$. As
indicated in \cite{ZhuBon1},
this definition guarantees that 
$d_{\mathrm H}(\lambda,
\lambda')$ is independent of the
metric of $S$ up to H\"older
equivalence, whereas it is
unclear whether the same
property holds for the Hausdorff
metric as closed subsets of $S$.
This subtlety is relevant only
when we consider metric
properties since, as proved in
\cite[Lemma~3.5]{CasBle},  the two
metrics define the same topology
on
$\mathcal L(S)$. 

A classical tool in
2--dimensional topology/geometry
is the notion of train track. A
\emph{train track} on  the
surface $S$
 is a graph $\Theta $ contained in
the interior of $S$ which 
consists of finitely many
vertices, also called 
\emph{switches}, and of finitely 
many  edges joining them such
that:
\begin{enumerate}
\item    The edges of $\Theta $
are  differentiable arcs whose
interiors are embedded  and
pairwise disjoint (the two end
points of an  edge may coincide).
\item    At each switch $s$ of
$\Theta
$,  the edges of $\Theta $ that contain $s$
are all  tangent to the same line $L_{s}$ in
the tangent  space $T_{s}S$ and, for each of
the two  directions of $L_{s}$, there is at
least one edge  which is tangent to that
direction. 
\item    Observe that the
complement 
$S-\Theta $ has a certain number of spikes,
each  leading to a switch $s$ and locally
delimited by  two edges that are tangent to
the same direction  at $s$; we require that no
component of $S-\Theta 
$ is a disc with 0, 1 or 2 spikes or an open 
annulus with no spike.
\end{enumerate}

A curve $c$
\emph{carried} by the train track
$\Theta$ is a differentiable
immersed curve  $c\co I
\rightarrow S$  whose image is
contained in
$\Theta$, where $I$ is an
interval in $\mathbb R$. The
geodesic lamination $\lambda$ is
\emph{weakly carried} by the
train track $\Theta$ if, for
every leaf $g$ of $\lambda$,
there is a curve $c$ carried by
$\Theta$ which is homotopic to
$g$ by a homotopy moving points
by a bounded amount. In this
case, the bi-infinite sequence
$\left\langle \dots, e_{-1}, e_0,
e_1, \dots, e_n, \dots
\right\rangle$ of the edges of
$\Theta$ that are crossed in this
order by the curve $c$ is the
\emph{edge path realized by} the
leaf $g$; it can be shown that
the curve
$c$ is uniquely determined by the
leaf
$g$, up to reparametrization, so
that the edge path realized by
$g$ is well-defined up to order
reversal. 

Let $\mathcal L(\Theta)$ be
the set of geodesic laminations
that are weakly carried by
$\Theta$. (This was denoted
by $\mathcal L^{\mathrm
w}(\Theta)$ in
\cite{ZhuBon1} where, unlike in
the current paper,  we had to
distinguish between ``strongly
carried'' and ``weakly carried''.)

We introduced two different
metrics on $\mathcal L(\Theta)$ in
\cite{ZhuBon1}. The first one is
defined over all of $\mathcal
L(S)$, and is just a variation
of the Hausdorff metric $d_{\mathrm
H}$. The distance
function $d_{\log}$ on
$\mathcal L(S)$ is defined by the
formula
\begin{equation*} d_{\log}
(\lambda,\lambda')=
\frac1 {\left\lvert \log 
\left(
\min\left\{ d_{\mathrm
H}(\lambda,\lambda'),
\frac14 \right\} 
\right)
\right\rvert}.
\end{equation*} 
In particular, 
$ d_{\log}
(\lambda,\lambda')=
{\left\lvert 1/\log 
 d_{\mathrm
H}(\lambda,\lambda')
\right\rvert}
$ when $\lambda$ and $\lambda'$ are
close enough from each other. The
$\min$ in the formula was only
introduced to make $d_{\log}$
satisfy the triangle inequality,
and is essentially cosmetic.  

The other metric is the
\emph{combinatorial distance}
between $\lambda$ and $\lambda'
\in \mathcal L(\Theta)$, defined
by 
\begin{equation*} 
\textstyle
d_\Theta(\lambda,
\lambda')=
\min \left\{ \frac1{r+1};
\lambda \text{ and } \lambda'
\text{ realize the same edge
paths of length } r
\right\},
\end{equation*}
where we say that an edge path is
realized by a geodesic lamination
when it is realized by one of its
leaves. This metric is actually
an ultrametric, in the sense that
it satisfies the stronger
triangle inequality 
$$d_\Theta
(\lambda, \lambda'') \leq \max
\left\{ d_\Theta
(\lambda, \lambda'), d_\Theta
(\lambda', \lambda'') \right\}.$$ 

The main interest of this
combinatorial distance is the
following fact, proved in
\cite{ZhuBon1}.
\begin{proposition}
\label{prop:d_thetaEqd_log}
For every train track $\Theta$ on
the surface $S$, the
combinatorial metric $d_\Theta$
is Lipschitz equivalent to the
restriction of the metric
$d_{\log}$ to  
$\mathcal L(\Theta)$.
\end{proposition}

The statement that $d_\Theta$ and
$d_{\log}$ are Lipschitz equivalent
means that there exists constants
$c_1$ and $c_2>0$ such that 
\begin{equation*}
c_1\, d_{\Theta}(\lambda, \lambda')
\leq  d_{\log}(\lambda, \lambda')
\leq c_2\, d_{\Theta}(\lambda,
\lambda')
\end{equation*}
for every $\lambda$, $\lambda'\in
\mathcal L(\Theta)$. In particular,
the two metrics define the same
topology, and have the same
Hausdorff dimension. 

For future reference,
we note the following property,
whose proof can be found in
\cite[Chapter~1]{Bon10} (and also
easily follows from
Proposition~\ref
{prop:d_thetaEqd_log}).

\begin{proposition}
\label{prop:L(Theta)compact}
The space $\mathcal L(\Theta)$ is
compact.
\end{proposition}

\section{Distance estimates on the
once-punctured torus}
\label{sect:EstimatesPunctTorus}

In this section, we focus
attention on the case where the
surface $S$ is the once-punctured
torus, which we will here denote by
$T$.  As indicated above, there is
a very convenient classification of
simple closed  geodesics or,
equivalently, isotopy classes of
simple closed curves,  on $T$; see
for instance
\cite{PenHar}. 

Let $\mathcal S(T)\subset
\mathcal L(T)$ denote the set of
all simple closed geodesics which
are contained in the interior of
$T$. In the plane
$\mathbb R^2$, consider the
lattice $\mathbb Z^2$. The
quotient of $\mathbb R^2 - \mathbb
Z^2$ under the group $\mathbb
Z^2$ acting by translations is
diffeomorphic to the interior of
$T$. Fix such an identification
$\mathrm {int}(T) \cong\left(
\mathbb R^2-
\mathbb Z^2\right) /\mathbb Z^2$. 
 Then every straight line in
$\mathbb R^2$ which has rational
slope and avoids $\mathbb Z^2$
projects to a simple closed curve
in
$\mathrm {int}(T)=\left( \mathbb
R^2-
\mathbb Z^2\right) /\mathbb Z^2$,
which itself is isotopic to a
unique simple closed geodesic of
$\mathcal S(T)$. This element of
$\mathcal S(T)$ depends only on
the slope of the line, and this
construction induces a bijection
$\mathcal S(T)\cong
\mathbb Q\cup \{\infty\}$. 

By definition, the element of 
$\mathbb Q\cup \{\infty\}$ thus
associated to $\lambda \in
\mathcal S(T)$ is the
\emph{slope} of $\lambda$. 

\begin{figure}[ht]
\centerline{
\includegraphics[scale=0.9]{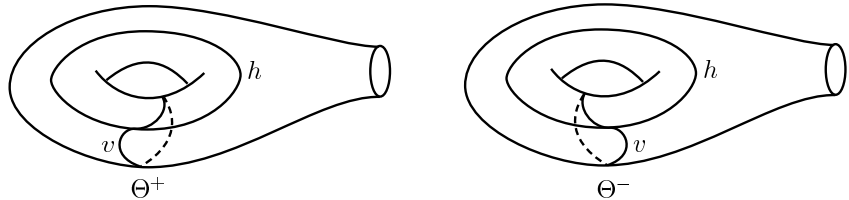}
}
\caption{The train tracks
$\Theta^+$ and $\Theta^-$ on the once-punctured torus $T$}
\label{pict:TrainTracksPunctTorus}
\end{figure}

\begin{figure}[ht]
\centerline{
\includegraphics{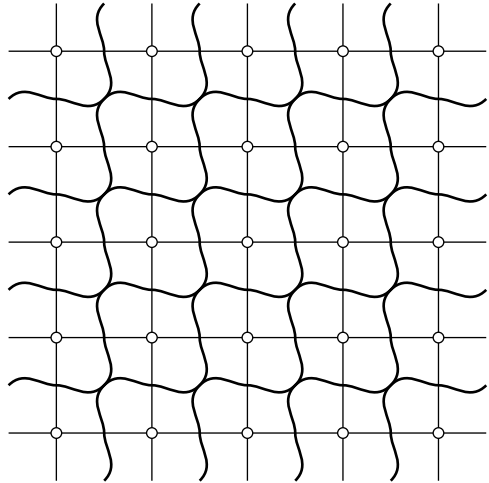}}
\caption{The preimage of
$\Theta^+$ in $\mathbb R^2-
\mathbb Z^2$}
\label{pict:R2-Z2}
\end{figure}

From this description, one
concludes that every simple closed
geodesic $\lambda \in
\mathcal S(T)$ is weakly carried
by one of the two train tracks
$\Theta^+$ and $\Theta^-$ represented
on Figure~\ref
{pict:TrainTracksPunctTorus}.
These two train tracks each
consist of two edges
$h$ and $v$ meeting at one single
switch. The identification
$\mathrm {int}(T) \cong\left(
\mathbb R^2-
\mathbb Z^2\right) /\mathbb Z^2$
can be chosen so that the preimage
of $\Theta^+$ in $\mathbb R^2 -
\mathbb Z^2$ is the one described
in Figure~\ref
{pict:R2-Z2}. In particular,
the preimage of the edge $h$ is a
family of `horizontal' curves,
each properly isotopic to a
horizontal line in $\mathbb R^2 -
\mathbb Z^2$, and the
preimage of the edge $v$  is a
family of `vertical' curves.
Similarly, the preimage of
$\Theta^-$ is obtained from that
of $\Theta^+$ by reflection
across the $x$--axis.

The simple closed
geodesic 
$\lambda
\in
\mathcal S(T)$ is weakly carried
by $\Theta^+$ (respectively $\Theta^-$)
exactly when its slope $\frac pq
\in
\mathbb Q\cup
\{\infty\}$ is non-negative
(respectively non-positive), by
consideration of a line of slope
$\frac pq$ in $\mathbb R^2 -
\mathbb Z^2$ and of its
translates under the action of
$\mathbb Z^2$. In this case, it is
tracked by a simple closed curve
$c$ carried by
$\Theta^+$ (respectively $\Theta^-$) which
crosses
$\left\lvert p\right\rvert$ times
the edge
$v$ and
$ q$ times the edge
$h$. We are here requiring the
integers $p$ and $q$ to be coprime
with
$q\geq 0$, and the slope $\infty =
\frac10= \frac{-1}0$ is considered
to be both non-negative and
non-positive. We will use the
same convention for slopes
throughout the paper. 

The following result, which
computes the combinatorial
distance between two simple closed
geodesics in terms of their
slopes, is the key to our
analysis of $\mathcal L_0(T)$. 

\begin{proposition}
\label{prop:d_phiAndSlopes} Let
the simple closed geodesics
$\lambda$, $\lambda'\in\mathcal
S(T)$ have slopes
$\frac pq$, $\frac {p'}{q'} \in
\mathbb Q\cup \{\infty\}$ with
$0\leq \frac pq  <\frac
{p'}{q'}\leq
\infty$. Then
\begin{equation*} 
\textstyle
d_{\Theta^+}
(\lambda, \lambda') = 
\max \left\{ \frac1{p''+q''};
 \frac pq \leq\frac {p''}{q''}
\leq\frac {p'}{q'}
\right\}.
\end{equation*} 
\end{proposition}

\begin{proof} For this, we first
have to understand the edge paths
realized by a simple closed
geodesic $\lambda\in\mathcal
S(T)$ in terms of its slope
$\frac pq$. 

Let $L$ be a line in $\mathbb
R^2$ of slope $\frac pq$ which
avoids the lattice $\mathbb Z^2$.
Look at its intersection points
with the grid $\mathbb Z\times
\mathbb R \,\cup\, \mathbb R
\times
\mathbb Z$, and label them as
$$\dots, x_{-1}, x_0, x_1,
\dots, x_i, \dots$$
in this
order along $L$. This defines a
periodic bi-infinite edge path
$$\left\langle \dots,  e_{-1},
e_0, e_1, \dots, e_i,
\dots \right\rangle$$ in $\Theta^+$,
where $e_i$ is equal to
the edge $h$ if the point $x_i$ is
in a vertical line
$\{n\}\times\mathbb R$ of the
grid, and $e_i=v$ if $x_i$ is
in a horizontal line $\mathbb
R\times \{n\}$. By consideration
of Figures~\ref
{pict:TrainTracksPunctTorus} and
\ref {pict:R2-Z2}, it is then
immediate that a (finite) edge
path is realized by
$\lambda$ if and only if it is
contained in this bi-infinite
edge path $\left\langle \dots, 
e_{-1}, e_0, e_1, \dots, e_i,
\dots \right\rangle$.

The main step in the proof of
Proposition~\ref
{prop:d_phiAndSlopes}
is the following special case.

\begin{lemma}
\label{thm:d_phiAndSlope1} If
$\lambda$, $\lambda'\in\mathcal
S(T)$ have finite positive slopes
$\frac pq$, $\frac {p'}{q'} \in
\mathbb Q
\cap\left]0,\infty\right[$ such
that $pq'-p'q=\pm 1$, then
$ d_{\Theta^+} (\lambda, \lambda')
= 
\max \left\{ \frac1{p+q},
\frac1{p'+q'} \right\}$.
\end{lemma}

\begin{proof}[Proof of
Lemma~\ref{thm:d_phiAndSlope1}]
Let $L$ and $L'$ be lines of
respective  slopes $\frac pq$ and
$\frac{p'}{q'}$ in $\mathbb R^2$,
avoiding the lattice $\mathbb
Z^2$. Let $c$ and $c'$ be the
projections of
$L$ and $L'$ to
$\mathrm{int}(T)\cong \left(
\mathbb R^2 -\mathbb Z^2 \right)/
\mathbb Z^2$. For suitable
orientations, the algebraic
intersection number of
$c$ and $c'$ is equal to
$pq'-p'q=\pm 1$. Since all
intersection points have the same
sign (depending on slopes and
orientations), we conclude that
$c$ and $c'$ meet in exactly one
point.

Let $A$ be the surface obtained
by splitting $\mathrm{int}(T)$
along the curve $c$.
Topologically, $A$ is a closed
annulus minus one point. Since
$c$ and $c'$ transversely meet in
one point, $c'$ gives in $A$ an
arc $c'_1$ going from one
component of $\partial A$ to the
other.

Similarly, the grid $\mathbb
Z\times
\mathbb R \,\cup\, \mathbb R
\times
\mathbb Z$ projects to a family
of arcs in $A$. Most of these
arcs go from one boundary
component of $A$ to the other.
However, exactly four of these
arcs go from $\partial A$ to the
puncture. We will call the union
of these four arcs the
\emph{cross} of $A$. As one goes
around the puncture, the arcs of
the cross are alternately
horizontal and vertical. Also,
the cross divides
$A$ into one hexagon and two
triangles
$\Delta$ and $\Delta'$. See
Figure~\ref{pict:SplitPunctTorus}
for the case where $\frac
pq=\frac 35$ and
$\frac{p'}{q'}=\frac 23$.

\begin{figure}[ht]
\centerline{
\includegraphics
{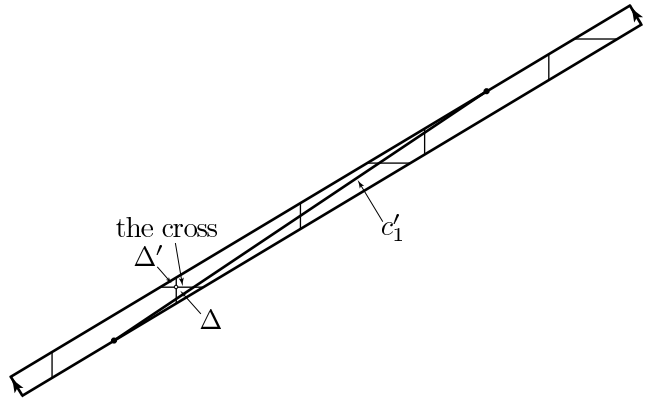}}
\caption{Splitting the once-punctured torus $T$ along the
image of a straight line
(glue the two short sides of
the rectangle)}
\label{pict:SplitPunctTorus}
\end{figure}

Set $r=\min\left\{p+q-1, p'+q'-1
\right\} \geq 0$. We want to show
that every edge path $\gamma=
\left\langle e_1, e_2, \dots,
e_{r}\right \rangle$ of length
$r$ in
$\Theta^+$ which is realized by
$c'$ is also realized by $c$. 

Given such an edge path $\gamma$, 
there exists an arc $a$ immersed
in
$c'$ which cuts the image of the
grid $\mathbb Z\times
\mathbb R \,\cup\, \mathbb R
\times
\mathbb Z$ at the points $x_1$,
$x_2$, \dots, $x_r$ in this
order, and such that $x_i$ is in
the image of a vertical line
$\{n\}\times\mathbb R$ if
$e_{i}=h$ and in the image of a
horizontal line
$\mathbb R\times\{n\}$ if
$e_{i}=v$. Since $c'$ crosses
the image of the grid in $p'+q'>r$
points, the arc $a'$ is actually
embedded in $c'$.

We had a degree of freedom in
choosing the closed curve $c$,
since it only needs to be the
projection of a line
$L$ with slope $\frac pq$. We can
choose this line $L$ so that $c$
contains the starting point of
the arc $a'$ (we may need to
slightly shorten $a'$ for this, in
order to make sure that $L
\subset \mathbb R^2$ avoids the
lattice
$\mathbb Z^2$). 

The arc $a'$ now projects to an
arc
$a_1'$ embedded in the arc
$c'_1\subset A$ traced by $c'$,
such that the starting point of
$a_1'$ is on the boundary of
$A$.  

Note that each boundary
component of $A$ crosses the
image of the grid $\mathbb Z\times
\mathbb R \,\cup\, \mathbb R
\times
\mathbb Z$ in $p+q>r$ points.
Since $a_1'$ cuts the image of
this grid in $r$ points, we
conclude that $a_1'$ ``turns less
than once'' around $A$, in the
sense that it cuts each arc of
the image of the grid in at most
one point. Similarly, if $\Delta$
and $\Delta'$ are the two
triangles delimited in $A$ by the
cross of the grid, $a_1'$ can
meet the union $\Delta\cup
\Delta'$ in at most one single
arc. It follows that there exists an
arc in $\partial A$ which cuts
exactly the same components of the
image of the grid as $a_1'$. 
This arc $a_1$
shows that the edge path $\gamma$
is also realized by
$c$, and therefore by $\lambda$. 

We conclude that every edge path
of length $r$ which is realized
by $\lambda'$ is realized by
$\lambda$. Exchanging the r\^oles
of
$\lambda$ and
$\lambda'$, every edge path of
length
$r$ which is realized by
$\lambda$ is also realized by
$\lambda'$. Consequently, $
d_{\Theta^+} (\lambda, \lambda')
\leq 
\frac 1{r+1}$.

To show that $ d_{\Theta^+}
(\lambda, \lambda') = \frac
1{r+1}$, we need to find an edge
path of length $r+1$ which is
realized by one of
$\lambda$, $\lambda'$ and not by
the other one. Without loss of
generality, we can assume that
$p+q\leq p'+q'$, so that
$r+1=p+q$. 

Consider $c$, $c'$, $A$ and
$c_1'$ as above. The curves $c$
and $c'$ cross the image of the
grid in $p+q$ and $p'+q'$ points,
respectively. By our hypothesis
that $p+q\leq p'+q'$, it follows
that $c'_1$ turns at least once
around the annulus $A$, and
therefore meets at least one of
the two triangles $\Delta$
delimited by the cross. By moving
the line
$L\subset\mathbb R^2$ projecting
to $c$ (while fixing $c'$ and the
corresponding line $L'$), we can
arrange that $\Delta\cap c_1'$
consists of a single arc, and
contains the initial point of
$c_1'$. Let $a_2'$ be an arc in
$c_1'$ which starts at this
initial point and crosses exactly
$p+q$ points of the grid. Note
that this is possible because
$p+q\leq p'+q'$. Because each of
the two components of $\partial A$
meets the grid in $p+q$ points,
the ending point of $a_2'$ is
contained in the other triangle
$\Delta'\not=\Delta$ delimited by
the cross. Consider the edge path
$\gamma'$ of length $p+q$
described by
$a'_2$.

By construction, the edge path
$\gamma'$ is realized by
$\lambda'$. We claim that it is
not realized by $\lambda$.
Indeed, let $\gamma$ be the edge path
described by an arc $a_2$ in
$\partial A$ which goes once
around the component of $\partial
A$ that contains the starting
point of $c_2'$, and starts and
ends at this point. By
construction, the edge path
$\gamma'$ is obtained from
$\gamma$ by switching the last
edge, either from $v$ to $h$, or
from
$h$ to $v$. In particular, the
edge paths $\gamma$ and $\gamma'$
contain different numbers of edges
$v$. However, because $c$ cuts
the grid in exactly $p+q$ points,
every edge path of length
$p+q$ which is realized by
$\lambda$ must contain the edge
$v$ exactly $p$ times. It follows
that $\gamma'$ is not realized by
$\lambda$.

We consequently found an edge
path $\gamma'$ of length $p+q=r+1$
which is realized by $\lambda'$
but not by $\lambda$. This proves
that $ d_{\Theta^+} (\lambda,
\lambda')
\geq 
\frac 1{r+1}$, and therefore that
$ d_{\Theta^+} (\lambda, \lambda')
=\frac 1{r+1}$. Since
$r=\min\left\{p+q-1, p'+q'-1
\right\}$, this concludes the
proof of
Lemma~\ref{thm:d_phiAndSlope1}.
\end{proof}

\begin{remark}
\label{rmk:d_phiAndSlope1} For
future reference, note that we
actually proved the following
property: Under the hypotheses of
Lemma~\ref{thm:d_phiAndSlope1}
and if $r=p+q-1\leq p'+q'-1$,
then $\lambda$ and $\lambda'$
realize exactly the same edge
paths of length $r$, and there
exists an edge path of length
$r+1$ which is realized by
$\lambda'$ and not by $\lambda$.
\end{remark}

\begin{lemma} 
\label
{thm:ConnectedSlopesForPath} If
$\lambda$,
$\lambda'$, $\lambda''\in\mathcal
S(T)$ have slopes
$\frac pq$, $\frac {p'}{q'}$,
$\frac {p''}{q''} \in
\mathbb Q$ with 
$0< \frac pq \leq\frac {p''}{q''}
\leq\frac {p'}{q'}<
\infty$, then every edge path
$\gamma$ in $\Theta^+$ which is
realized by both $\lambda$ and
$\lambda'$ is also realized by
$\lambda''$. 
\end{lemma}
\begin{proof}[Proof of Lemma~\ref
{thm:ConnectedSlopesForPath}]
 Let
$L$ and
$L'$ be lines of respective
slopes $\frac pq$ and
$\frac{p'}{q'}$ in $\mathbb
R^2-\mathbb Z^2$. Since $\lambda$
realizes $\gamma=\left\langle  e_1, e_2,
\dots, e_n
\right\rangle$, there is an
arc $a\subset L$ which meets the
grid $\mathbb Z\times
\mathbb R \,\cup\, \mathbb R
\times
\mathbb Z$ at the points $x_1$,
$x_2$, \dots, $x_n$ in this
order, and so that the point
$x_i$ is in a vertical line of
the grid when $e_i=h$ and in a
horizontal line when
$e_i=v$.  Since $\lambda'$
also realizes $\gamma'$, there is
a similar arc $a'\subset L'$ which
meets the grid  at points
$x_1'$,
$x_2'$, \dots, $x_n'$.

Applying to $L$ and $L'$ elements
of the translation group $\mathbb
Z^2$ if necessary, we can assume
without loss of generality that
the starting points of $a$ and
$a'$ are both in the square
$\left]0,1\right[\times
\left]0,1\right[$. Then, because
the slopes are both positive,  the
fact that the arcs
$a$ and
$a'$ cut the grid according to
the same vertical/horizontal
pattern implies that each $x_i$
is in the same line segment
component $I_i$ of 
$\left( \mathbb Z\times
\mathbb R \,\cup\, \mathbb R
\times
\mathbb Z \right) -\mathbb Z^2$
as $x_i'$. 

The set of lines which cut these
line segments $I_i$ is connected.
Therefore, one of them must have
slope
$\frac{p''}{q''}$. By a small
perturbation, we can arrange that
this line $L''$ of slope 
$\frac{p''}{q''}$ is also disjoint
from the lattice $\mathbb Z^2$.
The fact that  $L''$ cuts $I_1$,
$I_2$, \dots, $I_n$ in this order
then shows that the corresponding
simple closed geodesic
$\lambda''$ realizes the edge
path $\gamma$. 
\end{proof}

We can now conclude the proof of
Proposition~\ref
{prop:d_phiAndSlopes}. 
Temporarily setting aside the
slopes $0$ and $\infty$, let the
simple closed geodesics
$\lambda$, $\lambda'\in\mathcal
S(T)$ have slopes
$\frac pq$, $\frac {p'}{q'} \in
\mathbb Q$ with
$0< \frac pq  <\frac {p'}{q'}<
\infty$. Then, by elementary
number theory (see for instance
\cite[Section~3.1]{HarWri}), there is a
finite sequence of slopes 
$$\frac {p}{q} =
\frac {p_0}{q_0} <
\frac {p_1}{q_1}<
\dots<
\frac {p_n}{q_n} =
\frac {p'}{q'}$$ such that
$p_iq_{i-1}-p_{i-1}q_i=1$ for
every $i$. Let
$\lambda_i
\in
\mathcal S(T)$ be the simple
closed geodesic with slope $\frac
{p_i}{q_i}$

By the ultrametric property and by
Lemma~\ref{thm:d_phiAndSlope1}, 
\begin{equation*} 
\textstyle
d_{\Theta^+}
(\lambda, \lambda')
\leq  \max \left\{ d_{\Theta^+}
(\lambda_{i-1},
\lambda_i);\ i=1, \dots, n
\right\}=\frac1{r+1}
\end{equation*}  if $r=
\inf\left\{ p_i +q_i-1;\ i=0,
\dots, n
\right\}$. We want to prove that
this inequality is actually an
equality, namely that there is an
edge path of length $r+1$ which
is realized by one of
$\lambda$, $\lambda'$ and not by
the other. 

First consider the case where
$p+q-1>r$, and examine the first
$i$ such that $p_i+q_i-1=r$. By
Lemma~\ref{thm:d_phiAndSlope1} and
Remark~\ref {rmk:d_phiAndSlope1},
there is an edge path $\gamma$ of
length $r+1$ which is realized
by $\lambda$ and
$\lambda_{i-1}$, and not by
$\lambda_i$. Since $\frac pq
<\frac{p_i}{q_i} <\frac{p'}{q'}$,
Lemma~\ref
{thm:ConnectedSlopesForPath}
shows that $\gamma$ cannot be
realized by $\lambda'$, which
proves that 
$d_{\Theta^+} (\lambda, \lambda')
=\frac1{r+1}$.

When $p'+q'-1>r$, the same
argument provides an edge path
which is realized by $\lambda'$
and not by $\lambda$, again
showing that $d_{\Theta^+} (\lambda,
\lambda') =\frac1{r+1}$ in this
case.

Finally, consider the case where
$p+q-1=p'+q'-1=r$. Let $\gamma$
be any edge path of length $r+1$
which is realized by $\lambda$.
Note that
$\gamma$ goes exactly once around
$\lambda$. We conclude that
$\gamma$ contains exactly $p$
times the edge $v$, and $q$ times
the edge $h$. Similarly, any edge
path $\gamma'$ of length $r+1$
which is realized by $\lambda'$
must contain $p'$ times the edge
$v$, and $q'$ times the edge $h$.
Since $\frac pq\not =
\frac{p'}{q'}$, we conclude that
such a $\gamma'$ cannot be
realized by $\lambda$. Therefore, 
$d_{\Theta^+} (\lambda, \lambda')
=\frac1{r+1}$ again in this case.

This proves that 
\begin{equation*}
\begin{split}  
 d_{\Theta^+}
(\lambda,
\lambda')
&=\textstyle\frac1{r+1}\\ 
&\textstyle =
\max\left\{\frac1 {p_i
+q_i};\ i=0, \dots, n
\right\}\\ 
&\textstyle = \max\left\{\frac1
{p'' +q''};\ 
\frac pq
\leq\frac{p''}{q''}
\leq\frac{p'}{q'}
\right\}
\end{split}
\end{equation*} where the last
equality comes from the
elementary property that
$p''\geq p_{i-1}+p_i$ and 
$q''\geq q_{i-1}+q_i$ whenever 
$\frac{p_{i-1}}{q_{i-1}}
<\frac{p''}{q''}<
\frac{p_i}{q_i}$ (hint:
$p_iq_{i-1}-p_{i-1}q_i=1$).

This concludes the proof of
Proposition~\ref
{prop:d_phiAndSlopes} in the case
where $0<\frac pq<\frac{p'}{q'}
<\infty$.

	When, $\frac pq=\frac01$, note
that $\lambda$ never crosses the
edge $v$, but that $\lambda'$
does. This provides an edge path
of length 1 which is realized by
$\lambda'$ and not by $\lambda$.
Therefore, 
$d_{\Theta^+} (\lambda,
\lambda') =1 = \max\left\{\frac1
{p'' +q''};\ 
\frac 01
\leq\frac{p''}{q''}
\leq\frac{p'}{q'}
\right\}$ in this case as well.
The case where $\frac
{p'}{q'}=\infty =
\frac10$ is similar.
\end{proof}

\begin{corollary} 
\label{cor:d_logAndSlopes}
The slope map
$\mathcal S(T) \rightarrow
\mathbb Q\cup \{ \infty \}$ sends the
metric $d_{\log}$ to a metric
which is Lipschitz equivalent to
the metric $d$ on $\mathbb Q\cup
\{\infty\}$ defined by
\begin{equation*} 
\textstyle 
d\bigl(
\frac pq,
\frac{p'}{q'}\bigr) =\max
\bigl\{  \frac 1{\left\vert
p'' \right\vert +\lvert q''\rvert };
\frac pq
\leq
\frac{p''}{q''} \leq
\frac{p'}{q'} \bigr \}
\end{equation*}
for $\frac pq <\frac{p'}{q'}$. 
\end{corollary}
\begin{proof}
Propositions~\ref
{prop:d_thetaEqd_log} and \ref
{prop:d_phiAndSlopes} prove this
property for the restrictions of
$d_{\log}$ to
$\mathcal L(\Theta^+) \cap \mathcal
S(T)$ and
$\mathcal L(\Theta^-) \cap \mathcal
S(T)$. It
therefore suffices to show that
there is a positive lower bound
for the distances
$d_{\log}(\lambda,
\lambda')$ as $\lambda$,
$\lambda'$ range over all simple
closed geodesics such that
$\lambda$ has finite negative
slope and $\lambda'$ has finite
positive slope; indeed the
$d$--distance between the slopes
of such $\lambda$ and $\lambda'$
is equal to 1. 

We could prove this
geometrically, but we will
instead use Proposition~\ref
{prop:d_phiAndSlopes} and the
fact that the Lipschitz
equivalence class of $d_{\log}$
is invariant under
diffeomorphisms of $T$. Recall
that every diffeomorphism of $T$
acts on the slope set $\mathbb
Q\cup\{\infty\}$ by linear fractional
maps
$x\mapsto \frac{ax+b}{cx+d}$,
with $a$, $b$, $c$, $d\in \mathbb
Z$ and $ad-bc=\pm 1$, and that
every such linear fractional map
is realized by a diffeomorphism of
$T$. 

First consider the case where the
slope $\frac pq$ of $\lambda$ is
in the interval $\left[-1,
0\right[$. Let $\varphi_1$ be a
diffeomorphism of $T$ whose
action on the slopes is given by
$x\mapsto x+1=\frac {x+1}{0x+1}$.
Now $\varphi_1(\lambda)$ and
$\varphi_1(\lambda')$ both have
non-negative slopes, which are on
different sides of the number 1.
It follows from
Propositions~\ref
{prop:d_thetaEqd_log} and
\ref{prop:d_phiAndSlopes} that
$d_{\log}(\varphi_1(\lambda),
\varphi_1(\lambda'))\geq c_0$ for
some constant $c_0>0$. Since
$\varphi_1$ does not change the
Lipschitz class of $d_{\log}$, it
follows that there exists a
constant $c_1>0$ such that 
$d_{\log}(\lambda,
\lambda')\geq c_1
d_{\log}(\varphi_1(\lambda),
\varphi_1(\lambda'))$. Therefore,
$d_{\log}(\lambda,
\lambda')\geq c_1c_0$. 

Similarly, when $\frac pq$ is in
the interval $\left]\infty,
-1\right]$, consider the
diffeomorphism $\varphi_2$ of $T$
whose action on the slopes is
given by $x\mapsto \frac x{x+1}$.
The same argument as above gives
$d_{\log}(\lambda,
\lambda')\geq c_2c_0$. 
\end{proof}

\section{Chain-recurrent geodesic
laminations on the once-punctured
torus}
\label{sect:ChainRec}

A geodesic lamination $\lambda
\in \mathcal L(S)$ is
\emph{chain-recurrent} if it is in
the closure of the set of all
multicurves (consisting of finitely
many simple closed geodesics) in
$S$. See for instance
\cite[Chapter~1]{Bon10} for an
equivalent definition of
chain-recurrent geodesic laminations
which better explains  the
terminology.

When the surface $S$ is the
once-punctured torus $T$, a
multicurve is, either a simple
closed geodesic in the interior
fo $T$ (namely an element of
$\mathcal S(T)$), or the union of
$\partial T$ and of an element of
$\mathcal S(T)$, or just
$\partial T$. As a consequence, a
chain-recurrent geodesic
lamination in the interior of $T$
is a limit of simple closed
geodesics.

Let $\mathcal L_0^{\mathrm
{cr}}(T)$ denote the set of
chain-recurrent geodesic laminations
that are contained in the interior
of $T$. By the above remarks,
$\mathcal L_0^{\mathrm {cr}}(T)$ is
also the closure in
$\mathcal L_0(T)$ of the set
$\mathcal S(T)$ of all simple
closed geodesics. 

The space $\mathcal L(S)$ is
compact; see for instance 
\cite[Section~3]{CasBle}, \cite[Section~4.1]{CanEpsGre} or
\cite[Section~1.2]{Bon10}. Also, there is a
neighborhood $U$ of $\partial T$
such that every complete geodesic
meeting
$U$ must, either cross itself, or
be asymptotic to
$\partial T$, or be $\partial T$;
in particular, every geodesic
lamination which meets
$U$ must contain $\partial T$. It
follows that $\mathcal L_0(T)$ is
both open and closed in $\mathcal
L(T)$. As a consequence, $\mathcal
L_0^{\mathrm {cr}}(T)$ is compact.

We conclude that $\left( 
\mathcal L_0^{\mathrm
{cr}}(T), d_{\log} \right)$ is
the completion of $\left(
\mathcal S(T), d_{\log} \right)$.
By
Corollary~\ref
{cor:d_logAndSlopes}, 
$\left( 
\mathcal L_0^{\mathrm
{cr}}(T), d_{\log} \right)$ is
therefore Lipschitz equivalent to
the completion $\bigl ( \widehat
{\mathbb Q}, d \bigr)$ of 
$\bigl ( \mathbb Q \cup \{ \infty
\}, d
\bigr)$, where the metric $d$ is
defined by 
\begin{equation*} 
\textstyle 
d\bigl(
\frac pq,
\frac{p'}{q'}\bigr) =\max
\bigl\{  \frac 1{\lvert
p'' \rvert + \lvert q''\rvert};
\frac pq
\leq
\frac{p''}{q''} \leq
\frac{p'}{q'} \bigr \}
\end{equation*}
for $\frac pq <\frac{p'}{q'}$. 

This completion $\bigl ( \widehat
{\mathbb Q}, d \bigr)$ is studied
in detail in the Appendix. In
particular, Proposition~\ref
{prop:TopQhat}
determines its topology, and
Proposition~\ref
{prop:dim2} computes its
Hausdorff dimension and its
Hausdorff measure in this
dimension. These two results
prove:

\begin{theorem}
\label{thm:TopL_0cr}
The space $\mathcal L_0^{\mathrm
{cr}}(T)$ is
homeomorphic to the subspace $K\cup
L_1$ of the circle $\mathbb R \cup
\{\infty \}$ obtained by adding
to the standard
middle third Cantor set
$K \subset [0,1] \subset \mathbb
R$ a family $L_1$ of isolated points
consisting of exactly one point
 in each component of
$\mathbb R \cup
\{\infty \} - K$.
\end{theorem}

\begin{theorem}
\label{thm:DimL_0cr}
The metric space
$\left( 
\mathcal L_0^{\mathrm
{cr}}(T), d_{\log} \right)$ has
Hausdorff dimension $2$, and its
$2$--dimensional Hausdorff
measure is equal to $0$. 
\end{theorem}

\section{Dynamical properties of
geodesic laminations}
\label{sect:DynProp}

We collect in this section a few
general facts on geodesic
laminations which will be useful
to extend our analysis from
chain-recurrent geodesic laminations
to all geodesic laminations. 

 A
geodesic lamination
$\lambda$ is
\emph{recurrent} is every
half-leaf of
$\lambda$  comes back arbitrarily
close to its starting point, and
in the same direction. For
instance, a multicurve
(consisting of finitely many
disjoint simple closed geodesics)
is recurrent.

A geodesic lamination $\lambda$
cannot be recurrent if it contains
an \emph{infinite isolated leaf},
namely a leaf $g$ which is not
closed  and for which there
exists a small arc $k$ transverse
to $g$ such that $k\cap\lambda =
k\cap g$ consists of a single
point. 

\begin{proposition}
\label{thm:DynGeodLam} A geodesic
lamination $\lambda$ has finitely
many connected component. It can
be uniquely decomposed as the
union of a recurrent geodesic
lamination
$\lambda^{\mathrm r}$ and of
finitely many infinite isolated
leaves which spiral along
$\lambda^{\mathrm r}$. 
\end{proposition}
\begin{proof} See for instance
\cite[Theorem~4.2.8]{CanEpsGre}
or \cite[Chapter~1]{Bon10}.
\end{proof}

Here the statement that an
infinite leaf $g$ spirals along
$\lambda^{\mathrm r}$ means that
each half of $g$ is asymptotic to
a half-leaf contained in
$\lambda^{\mathrm r}$. 

Let a \emph{sink} in the geodesic
lamination $\lambda$ be an
oriented sublamination $\lambda_1
\subset \lambda$ such
that every half-leaf of
$\lambda-\lambda_1$ which spirals
along
$\lambda_1$ does so in the
direction of the orientation, and
such that there is at least one
such  half-leaf spiralling along
$\lambda_1$.

\begin{proposition}
{\rm\cite[Chapter~1]{Bon10}}\qua
\label{prop:NoSink}A geodesic lamination is
chain-recurrent if and only if it
contains no sink.
\end{proposition}

As a special case of
Proposition~\ref{prop:NoSink},
every recurrent geodesic
lamination is also chain-recurrent. 

In our analysis of the once-punctured torus and the 4--times-punctured
sphere, the following
lemma will be convenient to push
our arguments from chain-recurrent
geodesic laminations to all
geodesic laminations. We prove it
in full generality since it may
be of independent interest.

\begin{lemma} 
\label{lem:ConvNonRecurPart}
 There
exists constants $c_0$,
$r_0>0$, depending only
on the (negative) curvature of
the metric
$m$ on $S$, with the following
property. Let
$\lambda_1$ be a geodesic
lamination contained in the
geodesic lamination
$\lambda$ and containing the
recurrent part $\lambda^{\mathrm
r}$ of $\lambda$. Then any
geodesic lamination
$\lambda'_1$ with $d_{\mathrm
H}(\lambda_1,
\lambda_1')<r_0$ is
contained in a geodesic
lamination $\lambda'$ with 
$d_{\mathrm H}(\lambda,
\lambda')\leq c_0\,d_{\mathrm
H}(\lambda_1,
\lambda_1')$.
\end{lemma}

\begin{proof} We will explain how
to choose $c_0$ and
$r_0$ in the course
of the proof. Right now, assume
that
$r_0$ is given, and pick
$r$ with $
r /2 <d_{\mathrm
H}(\lambda_1,
\lambda_1')<r  \leq
r _0$.

We claim that there is a constant
$c_1>1$ such that, at each $x\in
\lambda \cap \lambda_1'$, 
the angle between the lines $T_x
\lambda$ and $T_x\lambda_1'$ is
bounded by $c_1 r$. 
Indeed,
since 
 $ d_{\mathrm H}(\lambda_1,
\lambda_1')<r $, there
is a point $y\in\lambda_1$ such
that $\bigl( y,
T_y\lambda_1\bigr)$ is at
distance less than $r $
from $\bigl( x,
T_{x}\lambda_1'\bigr)$ in the
projective tangent bundle
$PT(S)$. In particular, the
distance between the two points
$x$, $y\in\lambda$ is less than
$r $. In this situation,
a classical lemma (see
\cite[Corollary~2.5.2]{EpsMar} or
\cite[Appendix~B]{Bon10}) asserts
that, because the two leaves of
$\lambda$ passing through $x$ and
$y$ are disjoint or equal, 
there is a constant
$c_2$, depending only on the
curvature of the metric $m$, such
that the distance from $\bigl( x,
T_{x}\lambda\bigr)$ to $\bigl(
y, T_{y}\lambda\bigr) = \bigl(
y, T_{y}\lambda_1\bigr)$ in
$PT(S)$ is bounded by $c_2 d(x,
y)$, and therefore by
$c_2r $. Consequently,
the angle from $T_{x}\lambda$ to
$T_{x}\lambda_1'$ at $x$, namely the
distance from $\bigl( x,
T_{x}\lambda\bigr)$ to 
$\bigl( x,
T_{x}\lambda_1'\bigr)$ in
$PT(S)$, is bounded by
$(1+c_2)r=c_1r $. 

Since $\lambda_1$ contains the
recurrent part of $\lambda$,
Proposition~\ref{thm:DynGeodLam}
shows that $\lambda$ is the union
of $\lambda_1$ and of finitely
many infinite isolated leaves.
Let $\widehat\lambda_1'$ denote
the canonical lift of $\lambda_1'$
to the projective tangent bundle
$PT(S)$, consisting of those
$\bigl(x, T_x\lambda_1
\bigr) \in PT(S)$ where $x\in
\lambda_1$.
Let
$A$ consist of those points
$x$ in $\lambda-\lambda_1$ such
that $\bigl(x, T_x\lambda 
\bigr)$ is at distance greater than
$c_1 r $ from
$\widehat\lambda_1'$ in $PT(S)$,
where $c_1$ is the constant defined
above. The set
$A$ is disjoint from
$\lambda_1'$ by choice of $c_1$. 
Because 
$ d_{\mathrm H}(\lambda_1,
\lambda_1')<r < c_1r $ and
because the leaves of
$\lambda-\lambda_1$ spiral along
the recurrent part
$\lambda^{\mathrm r}\subset
\lambda_1$, the set $A$ stays
away from the ends of
$\lambda-\lambda_1$. As a
consequence, $A$ has only
finitely many components $a_1$,
$a_2$, \dots, $a_n$ whose length
is at least
$r $.

Let us focus attention on one of
these $a_i$, contained in an
infinite isolated leaf
$g_i$ of $\lambda-\lambda_1$. Let
$b_i$ be the component of
$g_i-\lambda_1'$ that contains
$a_i$. The open interval $b_i$
can have 0, 1 or 2 end points in
$g_i$ (corresponding to points
where $\lambda_1'$ transversely
cuts $g_i$). 

Let $x_i$ be an end point of
$b_i$. Then $x_i$ is contained in
a leaf $g_i'$ of $\lambda_1'$. 
We observed that the angle between
$g_i$ and
$g_i'$ at $x_i$ is bounded by
$c_1r$.  Let $k_i$
be the half-leaf of $\lambda_1'$
delimited by $x_i$ in $g_i'$
which makes an angle of at least
$\pi-c_1r $ with
$b_i$ at $x_i$; note that $k_i$
is uniquely determined if we
choose $r _0$ small enough that
$c_1r \leq c_1r _0
<\pi/2$. 

We now construct a family $h$ of
bi-infinite or closed piecewise
geodesics such that:
\begin{enumerate}
\item $h$ is the union of all the
arcs $b_i$ and of pieces of the
half-leaves $k_j$ considered
above.
\item The external angle of $h_i$
at each corner is at most
$c_1r$.
\item $h$ can be perturbed to a
family of disjoint simple curves
contained in the complement of
$\lambda_1'$.
\item One of the two geodesic
pieces meeting at each
corner of $h$ has length at least
$1$ (say).
\end{enumerate}
As a first approximation and if
we do not worry about the third
condition, we can just take $h$
to be the union of the arcs 
$b_i$ and of the half-leaves $k_j$
(of infinite length). However,
with respect to this third
condition, a problem arises when
one half-leaf
$k_i$ collides with another
$k_j$; more precisely when, as
one follows the half-leaf $k_i$
away from the end point $x_i$ of
$b_i$, one meets an end point
$x_j$ of another arc
$b_j$ (with possibly $b_j=b_i$)
such that $b_j$ is on the same
side of
$k_i$ as $b_i$ and such that the
half-leaf $k_j$ associated to
$x_j$ goes in the direction
opposite to $k_i$. In this
situation, remove from $h$ the
two half-leaves $k_i$ and $k_j$
and add the arc $c_{ij}$
connecting
$x_i$ to $x_j$ in $k_i$. Because
the leaves of $\lambda$
containing $b_i$ and $b_j$ do not
cross each other, the length of
$c_{ij}$ will be at least $1$ if
we choose $r_0$ so that $c_1r
\leq c_1r_0$ is small enough,
depending on the curvature of
the metric.

Iterating this process, one
eventually reaches an $h$
satisfying the required
conditions.

Consider a component $h_i$ of
$h$. By construction, $h_i$ is
piecewise geodesic, the external
angles at its corners are at most
$c_1r$, and  every other
straight piece of $h_i$ has
length at least $1$. 
A Jacobi field argument then
provides a constant $c_3$ such
that
$h_i$ can be deformed to a
geodesic
$h_i'$ by a homotopy which moves
points by a distance bounded by
$c_3r $. Actually, a
little more holds: if the
homotopy sends $x\in h_i$ to
$x\in h_i'$, then the distance
from $\bigl(x , T_x h_i\bigr)$ to 
$\bigl(x' , T_{x'} h_i'\bigr)$ in
$PT(S)$ is bounded by
$c_3r $.

Consider the
geodesics $h_i'$ thus associated
to the components $h_i$ of $h$. By
the Condition~(3)
imposed on $h$, the
$h_i$ are simple, two
$h_i$ and $h_j$ are
either disjoint or
equal, and each $h_i$ is either
disjoint from
$\lambda_1'$ or contained in it.  
Also, by construction of
$h$, each end of a geodesic
$h_i'$ which is not closed is
 asymptotic either to a leaf of
$\lambda_1'$ (containing a
half-leaf $k_j$) or to a leaf of
$\lambda^{\mathrm r}$ which is
disjoint from $\lambda_1'$
(and containing an infinite arc
$b_j$). 

Let $\lambda'$ be the union of
$\lambda_1'$ and of the closure
of the geodesics  $h_i'$
thus associated to the components
$a_i$ of length $\geq r$ of $A$. By
the above observations,
$\lambda'$ is a geodesic
lamination. 

 We want to prove that
$d_{\mathrm H}(\lambda, \lambda')
\leq c_0 r $ for some
constant
$c_0$. Let $\widehat \lambda$, 
$\widehat \lambda'$,
$\widehat \lambda_1$ and 
$\widehat \lambda_1'$ denote the
respective lifts of $\lambda$,
$\lambda'$, $\lambda_1$ and
$\lambda_1'$ to the projective
tangent bundle $PT(S)$. 

If $x'$ is a point of $\lambda'$,
either it belongs to
$\lambda_1'$, or it belongs to
one of the geodesics
$h_i'$, or it belongs to one of
the components of
$\lambda^{\mathrm r}$ (in the
closure of some $h_i'$). In the
first and last case, it is
immediate that the corresponding
point $\bigl(x' , T_{x'}
\lambda'\bigr)\in \widehat
\lambda'$ is at distance less than
$r $ from
$\widehat\lambda$ in $PT(S)$. If
$x'$ is in the geodesic $h_i'$,
then we saw that there is a point
$x\in h_i$ such that the distance
from $\bigl(x' , T_{x'}
\lambda'\bigr)=\bigl(x' , T_{x'}
h_i'\bigr)$ to $\bigl(x, T_{x}
h_i\bigr)$ is at most
$c_3r $. Then $\bigl(x,
T_{x} h_i\bigr)$ belongs to
$\widehat\lambda$ if $x$ is some
arc
$b_j$, and belongs to
$\widehat\lambda_1'$ otherwise.
Since  $ d_{\mathrm
H}(\widehat\lambda_1, \widehat
\lambda_1')<r $, we
conclude that 
$\bigl(x' , T_{x'}
\lambda'\bigr)$ is at distance at
most $(c_3+1)r $ from
$\widehat\lambda$ in this case.
This proves that
$\widehat\lambda'$ is contained
in the
$(c_3+1)r $--neighborhood
of $\widehat\lambda$.

Conversely, if $x$ is a point of
$\lambda$, either $\bigl(x,
T_x\lambda 
\bigr)$ is at distance at most 
$c_1 r $ from
$\widehat\lambda_1'\subset
\widehat\lambda'$, or $x$ belongs
to the subset $A$ of
$\lambda-\lambda_1$ introduced at
the beginning of this proof. If
$x$ belongs to one of the
components $a_i$ of $A$ used to
construct the leaves $h_i'$ of
$\lambda'$, then
$\bigl(x, T_x\lambda 
\bigr)$ is at distance less than
$c_3r $ from 
$\bigl(x, T_x h_i'
\bigr) \in \widehat\lambda'$ for
some $x\in h_i' \subset\lambda'$.
If $x$ belongs to a component
$a$ (of length $<r $) of
$A$ which is not one of the
$a_i$, then $x$ is at distance
less than $\frac 12 r $ from an
end point $y$ of $a$; in this
case 
$\bigl(x, T_x\lambda 
\bigr)$ is at distance less than
$\frac12 r $ from $\bigl(y,
T_y\lambda 
\bigr)$ by definition of the
metric of $PT(S)$, and $\bigl(y,
T_y\lambda 
\bigr)$ is at distance
$c_1 r $ from
$\widehat\lambda_1'\subset
\widehat\lambda'$ by definition
of $A$. We conclude that
$\bigl(x, T_x\lambda 
\bigr)$ is at distance at most
$\max\{ (c_1 + \frac12) r ,
c_3r \}$ from
$\widehat\lambda'$ in all cases.
Consequently, $\widehat\lambda$
is contained in the
$\max\{ (c_1 + \frac12) r ,
c_3r \}$--neighborhood
of $\widehat\lambda'$. 

This proves that, if we set 
 $c_0=2
\max\bigl\{ c_1 + \frac12 ,
c_3+1\bigr\}$, then
$d_{\mathrm H}(\lambda,
\lambda')=d_{\mathrm
H}(\widehat\lambda,
\widehat\lambda')\leq
c_0r /2 <c_0 d_{\mathrm
H}(\lambda_1,
\lambda_1')$ by choice of
$r$.  
\end{proof}

\section{The topology of geodesic
laminations of the once-punctured
torus}
\label{sect:TopGeodLamPuncTor}

Every recurrent geodesic
lamination admits a full support
transverse measure. The following
is a consequence of the fact that
there is a relatively simple
classification of measured
geodesic laminations on the once-punctured torus. See for instance
\cite{PenHar}, or \cite
{FatLauPoe} using the closely
related notion of measured
foliations.

\begin{proposition}
\label{prop:RecGeodLamPunctTor}
Every recurrent geodesic
lamination in the interior of the
once-punctured torus $T$ is
orientable, and admits a unique
transverse measure up to
multiplication by a positive real
number. This establishes a
correspondence between the set of
recurrent geodesic laminations in
the interior of $T$ and the set
of lines passing through the
origin in the homology space
$H_1(T; \mathbb R)$. 

When $\lambda$ corresponds to a
rational line, namely to a
line passing through non-zero
points of
$H_1(T; \mathbb Z)\subset H_1(T;
\mathbb R)$, the geodesic
lamination $\lambda$ is a simple
closed geodesic, and the
completion of its complement is a
once-punctured open annulus.
Otherwise, $\lambda$ has
uncountably many leaves and the
completion of its complement is a
once-punctured bigon, with two
infinite spikes. 
\end{proposition}

In this statement, the completion
of the complement $S-\lambda$ of
a geodesic lamination $\lambda$ in
a surface $S$ means its completion
for the path metric induced by the
metric of
$S$. It is always a
surface with geodesic boundary
and with finite area, possibly
with a finite number of infinite
spikes. See for instance
\cite[Section~4.2]{CanEpsGre} or
\cite[Chapter~1]{Bon10}. For instance,
when $\lambda$ corresponds to an
irrational line in $H_1(T;\mathbb
Z)$ in
Proposition~\ref
{prop:RecGeodLamPunctTor},
the completion of
$T-\lambda$ topologically is a
closed annulus minus two points
on one of its boundary components;
the boundary of this completion
consists of
$\partial T$ and of two geodesics
corresponding to infinite leaves of
$\lambda$ and whose ends are
separated by two infinite spikes. 

Fix an identification of the
interior of the once-punctured
torus with $\bigl( \mathbb R^2 -
\mathbb Z^2 \bigr) /\mathbb Z^2$.
This determines an identification
$H_1(T; \mathbb R) \cong \mathbb
R^2$, and a line in $H_1(T;
\mathbb R)$ is now determined by
its slope $s \in \mathbb
R\cup\{\infty\}$. Let $\lambda_s$
be the recurrent geodesic
lamination associated to the line
of slope $s$. 

The identification
$\mathrm{int}(T)\cong 
\bigl( \mathbb R^2 -
\mathbb Z^2 \bigr) /\mathbb Z^2$
also determines an orientation
for $T$. 

An immediate corollary of
Proposition~\ref
{prop:RecGeodLamPunctTor} is
that, if $\lambda$ is a geodesic
lamination in the interior of $T$
with recurrent part
$\lambda^{\mathrm r}$, the
completion of $T-\lambda^{\mathrm
r}$ contains only finitely many
simple geodesics (finite or
infinite). As a consequence, if
we are given the recurrent part
$\lambda^{\mathrm r}$, there are
only finitely many possibilities
for $\lambda$ and it is a simple
exercise to list all of them. 
Using
Proposition~\ref{prop:NoSink}, we
begin by enumerating the
possibilities for chain-recurrent
geodesic laminations.

\begin{proposition}
\label{prop:ChainRecPunctTor}
The chain-recurrent geodesic
laminations in the interior of the
once-punctured torus $T$ fall
into the following categories:
\begin{enumerate}
\item The recurrent geodesic
lamination
$\lambda_s$ with irrational
slope  $s \in \mathbb R-\mathbb
Q$.
\item The simple closed geodesic
$\lambda_s$ with rational slope
$s\in
\mathbb Q\cup \{\infty\}$.
\item The union $\lambda_s^+$ of
the simple closed geodesic 
$\lambda_s$ with slope $s\in
\mathbb Q\cup \{\infty\}$ and of
one infinite geodesic $g$ such
that, for an arbitrary
orientation of $\lambda_s$, one
end of $g$ spirals on the
right side of $g$ in the direction
of the orientation and the other
end spirals on the left side of
$g$ in the opposite
direction.
\item The union $\lambda_s^-$ of
the simple closed geodesic 
$\lambda_s$ with slope $s\in
\mathbb Q\cup \{\infty\}$ and of
one infinite geodesic $g$ such
that, for an arbitrary
orientation of $\lambda_s$, one
end of $g$ spirals on the
left side of $g$ in the direction
of the orientation and the other
end spirals on the right side of
$g$ in the opposite
direction.
\end{enumerate}
\end{proposition}

\begin{figure}[ht!]
\centerline{
\includegraphics[scale=0.7]{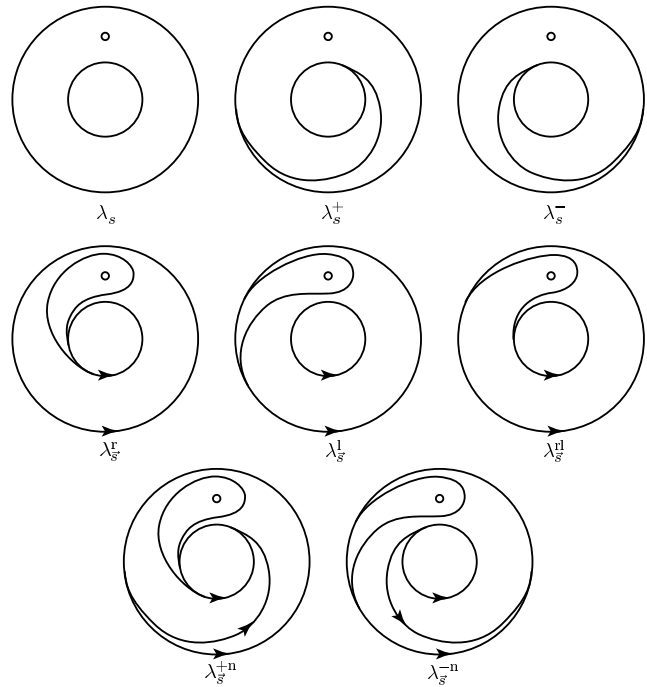}}
\caption{The geodesic laminations
containing the simple closed
geodesic
$\lambda_{s}
$, as
seen in the completion of
$T-\lambda_s$}
\label{pict:FiniteGeodLam}
\end{figure}

Note that, in Cases~3 and 4,
reversing the orientation  of
$\lambda_s$ exchanges left and
right, so that $\lambda_s^+$ and
$\lambda_s^-$ do not depend on
the choice of orientation for
$\lambda_s$. These two geodesic
laminations are illustrated in
Figure~\ref{pict:FiniteGeodLam}. 

A corollary of Proposition~\ref
{prop:ChainRecPunctTor} is that
every chain-recurrent geodesic
lamination in the interior of the
punctured torus is connected and
orientable.

In Theorem~\ref{thm:TopL_0cr}
(based on
Proposition~\ref{prop:TopQhat} in
the Appendix), we constructed a
homeomorphism $\varphi$ from the
space
$\mathcal L_0^{\mathrm{cr}}(T)$
of the chain-recurrent geodesic
laminations to the subspace
$K \cup L_1$ of $\mathbb R
\cup \{\infty\}$ union of the
standard middle third Cantor set
$K \subset [0,1]$ and of a family
$L_1$ of isolated points consisting
of the point
$\infty$ and of exactly one point
in each component of
$[0,1]-K$. We can revisit this
construction within the framework
of Proposition~\ref
{prop:ChainRecPunctTor}. 

\begin{proposition}
\label{prop:TopL_0Bis}
The homeomorphism $\varphi\co
\mathcal L_0^{\mathrm{cr}}(T)
\rightarrow K \cup L_1$
constructed in the proof of
Theorem~\ref{thm:TopL_0cr} is such
that:
\begin{enumerate}
\item [\rm(i)]The image under
$\varphi$ of the subset
$\mathcal S(T)\subset
\mathcal L_0^{\mathrm{cr}}(T)$
of simple closed curves is
exactly the set $L_1$ of isolated
points.
\item [\rm(ii)] If $s$, $t\in \mathbb
Q$ are such that $s<t$, then
$\varphi(\lambda_s)<
\varphi(\lambda_t)$ in $\mathbb
R$.
\item [\rm(iii)] If $I_s$ is the
component of $[0,1]-K$ containing
the image 
$\varphi\left(
\lambda_s \right)$ of the
simple closed geodesic
$\lambda_s$ of finite slope $s\in
\mathbb Q$, then, in the notation
of Proposition~\ref
{prop:ChainRecPunctTor}, the
left end point of
$I_s$ is
$\varphi\left(
\lambda_s^- \right)$ and its
right end
point is $\varphi\left(
\lambda_s^+ \right)$.
\item [\rm(iv)] $\varphi\left(
\lambda_\infty \right)=\infty$,
$\varphi\left(
\lambda_\infty^- \right)=1$ and
$\varphi\left(
\lambda_\infty^+ \right)=0$.
\end{enumerate}
\end{proposition}

\begin{proof} Properties~(i) and
(ii) are immediate from the
construction of $\varphi$ in
Theorem~\ref{thm:TopL_0cr} and
Proposition~\ref{prop:TopQhat}.
Note that Property~(i) is 
satisfied by an arbitrary
homeomorphism (since a
homeomorphism sends isolated
point to isolated point), but
that this is false for
Property~(ii). 

 A consequence of the
order-preserving condition of
Property~(ii) is that the left
end point of the interval $I_s$
is the limit of
$\varphi (\lambda_t)$ as $t$
tends to $s$ on the left. Since
$\lambda_t$ tends to
$\lambda_s^-$ as $t$ tends to $s$
on the left, it follows by
continuity of
$\varphi$ that the left end
point of $I_s$ is equal to
$\varphi(\lambda_s^-)$.
Similarly, the right end point
of $I_s$ is equal to the limit of 
$\varphi (\lambda_t)$ as $t$
tends to $s$ on the right, namely 
$\varphi(\lambda_s^+)$.

At $s=\infty$,
$\varphi(\lambda_\infty)=\infty$
by construction. By Property~(ii),
the point
$1$ is equal to the limit of
$\varphi(\lambda_t)$ as $t$ tends
to $\infty$ on the left (namely as
$t$ tends to $+\infty$). As $t$
tends to $\infty$ on the left,
$\lambda_t$ tends to
$\lambda_\infty^-$, and it
follows that
$\varphi(\lambda_\infty^-)=1$.
Similarly, $0$ is equal to the
limit of $\varphi(\lambda_t)$ as $t$ tends
to $\infty$ on the right (namely
as
$t$ tends to $-\infty$), and is
therefore equal to
$\varphi(\lambda_\infty^+)$. 
\end{proof}

We saw that a recurrent geodesic
lamination is connected
and orientable, and that
it is classified by its slope
$s\in
\mathbb R\cup \{\infty\}$. We can
interpret the slope as an element
of the space of unoriented lines
passing through the origin in
$\mathbb R^2$, namely as an
element of the projective plane
$\mathbb {RP}^1 = \mathbb R \cup
\{\infty\}$.

We will need to consider the
space of \emph{oriented}
recurrent geodesic laminations.
This clearly is a 2--fold cover of
the space of unoriented recurrent
geodesic laminations, and an
oriented recurrent geodesic
lamination is therefore
classified by its \emph{oriented
slope} $\vec s$,
defined as an element of the
space of oriented lines passing
through the origin in $\mathbb
R^2$, namely as an element of the
unit circle $\mathbb S^1$ in
$\mathbb R^2$. Let
$\lambda_{
\vec s}$ be the
oriented recurrent geodesic
lamination associated to the
oriented slope
$\vec s$. We similarly
define the oriented chain-recurrent
geodesic lamination
$\lambda_{\vec s}^+$
and 
$\lambda_{\vec
s}^-$ associated to the
irrational oriented slope
$\vec s$, union of the
oriented geodesic lamination
$\lambda_{\vec s}$
and of one additional geodesic
as in Cases~3 and 4 of
Proposition~\ref
{prop:ChainRecPunctTor}. 

We will say that the oriented
slope $\vec s \in
\mathbb S^1$ is
\emph{rational} when its
associated unoriented slope $s
\in\mathbb {RP}^1 = \mathbb R
\cup \{\infty\}$ is rational.

Again, a case-by-case analysis
combining Propositions~\ref
{prop:NoSink} and
\ref{prop:RecGeodLamPunctTor}
provides:
\begin{proposition}
\label{prop:NonChainRecPunctTor}
The geodesic laminations in the
interior of the punctured torus
$T$ which are not chain-recurrent
fall into the following
categories:
\begin{enumerate}

\item The union 
$\lambda_{\vec
s}^{\mathrm n}$ of
the oriented geodesic lamination
$\lambda_{\vec s}$,
with irrational oriented slope
$\vec s$,
and of one additional geodesic g
whose two ends converge to the
same spike of $T-\lambda_{
\vec s}$, in the
direction given by the
orientation.

\item The union 
$\lambda_{\vec
s}^{\mathrm r}$ of
the oriented simple closed
geodesic
$\lambda_{\vec s}$,
with rational oriented slope
$\vec s$,
and of one additional geodesic g
whose two ends spiral on the right
side of $\lambda_{\vec
s}$ in the direction given by the
orientation.

\item The union 
$\lambda_{\vec
s}^{\mathrm l}$ of
the oriented simple closed
geodesic
$\lambda_{\vec s}$,
with rational oriented slope
$\vec s$,
and of one additional geodesic g
whose two ends spiral on the left
side of $\lambda_{\vec
s}$ in the direction given by the
orientation.

\item The union 
$\lambda_{\vec
s}^{\mathrm {rl}}$ of
the oriented simple closed
geodesic
$\lambda_{\vec s}$,
with rational oriented slope
$\vec s$,
and of one additional geodesic g
whose two ends spiral around
$\lambda_{\vec s}$ in
the direction given by the
orientation, one on the right
side and one on the left side.

\item The union 
$\lambda_{\vec
s}^{+\mathrm n}$ of
the oriented chain-recurrent
geodesic lamination
$\lambda_{\vec s}^+$,
with rational oriented slope
$\vec s$,
and of one additional geodesic $g$
whose two ends converge to the
same spike of $T-\lambda_{
\vec s}^+$, in the
direction of the orientation.

\item The union 
$\lambda_{\vec
s}^{-\mathrm n}$ of
the oriented chain-recurrent
geodesic lamination
$\lambda_{\vec s}^-$,
with rational oriented slope
$\vec s$,
and of one additional geodesic $g$
whose two ends converge to the
same spike of $T-\lambda_{
\vec s}^-$, in the
direction given by the
orientation.
\end{enumerate}
\end{proposition}

\begin{figure}[ht!]
\centerline{
\includegraphics[scale=0.7]{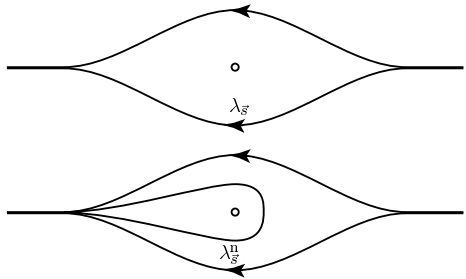}}

\caption{The geodesic laminations
$\lambda_{\vec
s}$ and
$\lambda_{\vec
s}^{\mathrm n}$ with $\vec s$
irrational, as seen in the
completion of
$T-\lambda_{\vec
s}$}
\label{pict:NonRecIrrat}
\end{figure}

Here the letters $\mathrm n$,
$\mathrm r$ and $\mathrm l$
respectively stand for
``non-chain-recurrent'', ``right''
and ``left''. Figure~\ref
{pict:NonRecIrrat} shows the
geodesic lamination 
$\lambda_{\vec
s}^{\mathrm n}$, for an
irrational oriented slope
$\vec s$. Figure~\ref
{pict:FiniteGeodLam} illustrates
the geodesic laminations
$\lambda_{\vec
s}^{\mathrm r}$,
$\lambda_{\vec
s}^{\mathrm l}$,
$\lambda_{\vec
s}^{\mathrm {rl}}$,
$\lambda_{\vec
s}^{+\mathrm n}$ and
$\lambda_{\vec
s}^{-\mathrm n}$ when the oriented
slope
$\vec s$ corresponds
to one orientation of the
rational slope
$s$. 

\section{The topology of
$\mathcal{L}_0(T)$\! for the
once-punctured torus}

We already determined the
topology of the space $\mathcal
L_0^{\mathrm {cr}}(T)$ of all
chain-recurrent geodesic laminations
in Theorem~\ref{thm:TopL_0cr}.  Let
$\mathcal L^{\mathrm n}_0(T)$
denote its complement $\mathcal
L_0(T)-\mathcal L_0^{\mathrm
{cr}}(T)$, namely the space  of all
non-chain-recurrent geodesic
laminations in the interior of
the once-punctured torus $T$. 

The following property is very
specific to the once-punctured
torus. For instance, we will see
in Section~\ref{sect:PunctSphere}
that it is false for the 4--times-punctured sphere.

\begin{lemma}
\label{lem:NonRecCompact}
The space 
 $\mathcal L^{\mathrm n}_0(T)$ is
compact. 
\end{lemma}
\begin{proof} By Proposition~\ref
{prop:NoSink} and by inspection
in  Proposition~\ref
{prop:NonChainRecPunctTor},
a geodesic lamination $\lambda$ on
the once-punctured torus which is
not chain-recurrent admits a
unique decomposition as the union
of an oriented chain-recurrent
 geodesic
lamination $\sigma(\lambda)$
(its unique sink) and of exactly
one infinite isolated leaf whose
two ends spiral along
$\sigma(\lambda)$ in the
direction of the orientation. The
chain-recurrent geodesic
lamination
$\sigma(\lambda)$ is weakly
carried by one of the two train
tracks $\Theta^+$ and $\Theta^-$
of
Figure~\ref
{pict:TrainTracksPunctTorus}.
By inspection, it follows that
$\lambda$ is weakly carried by
one of the four train tracks
$\Theta^\pm_\pm$  of Figure~\ref
{fig:FourNonRecTrainTracks},
unless $\lambda$ is of the form
$\lambda_{\vec s}^
{\mathrm{rl}}$ where the oriented
slope $\vec s$
corresponds to the unoriented
slopes 0 or $\infty$. The train
track $\Theta^\pm_\pm$  is made
up of
$\Theta^\pm$ and of one
additional edge $e$ going from one
``armpit'' of
$\Theta^\pm$ to itself; in
addition the edges of
$\Theta^\pm$ are oriented in such
a way that the orientations match
at the switch, and that the two
ends of the additional edge $e$
merge with $\Theta^\pm$ in the
direction of the orientation.
Note that the infinite isolated
leaf of $\lambda$ is tracked by a
curve carried by $\Theta^\pm_\pm$
which crosses the edge $e$
exactly once.

\begin{figure}[ht]
\centerline{
\includegraphics[scale=0.8]{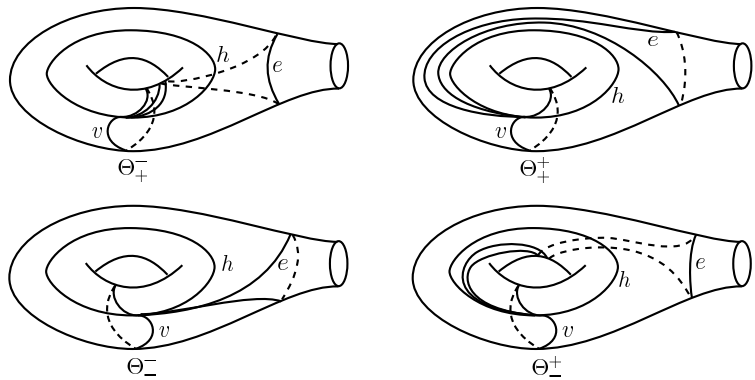}}
\caption{The train tracks
$\Theta^\pm_\pm$}
\label{fig:FourNonRecTrainTracks}
\end{figure}

Let $\mathcal
L^{\mathrm{n}}(\Theta^\pm_\pm)$
denote the space of
non-chain-recurrent geodesic
laminations that are weakly carried
by the train track
$\Theta^\pm_\pm$. We claim that
$\mathcal
L^{\mathrm{n}}(\Theta^\pm_\pm)$ is
equal to the set of those
$\lambda \in \mathcal
L(\Theta^\pm_\pm)$
(namely those $\lambda$ which are
weakly carried by
$\Theta^\pm_\pm$) which realize
the edge path $\left\langle e
\right\rangle$ consisting of the
single edge $e$. Indeed, a curve
which is carried by
$\Theta^\pm_\pm$ can cross $e$
at most once. If $\lambda \in \mathcal
L(\Theta^\pm_\pm)$
realizes $\left\langle e
\right\rangle$, it follows that
every leaf of $\lambda$ realizing $\left\langle
e
\right\rangle$ is isolated, and
consequently that those leaves
which do not realize $\left\langle
e
\right\rangle$ form a (closed)
non-empty sublamination
$\lambda^1 \subset \lambda$. Since
the train track
$\Theta^\pm \subset
\Theta^\pm_\pm$ is oriented, we
conclude that $\lambda_1$ is a
sink for $\lambda$, and therefore
that $\lambda$ is not
chain-recurrent by Proposition~\ref
{prop:NoSink}. 

A corollary of this observation
is that $\mathcal
L^{\mathrm{n}}(\Theta^\pm_\pm)$
is closed in $\mathcal
L(\Theta^\pm_\pm)$
(since the topology can be
defined with the metric
$d_{\Theta_\pm^\pm}$). Because
$\mathcal
L(\Theta^\pm_\pm)$
is compact by
Proposition~\ref
{prop:L(Theta)compact},
it follows that $\mathcal
L^{\mathrm{n}}(\Theta^\pm_\pm)$
is compact.

As a consequence, the space
\begin{equation*}
 \mathcal
L^{\mathrm{n}}(T)=
\mathcal
L^{\mathrm{n}}(\Theta^+_+)
\cup
\mathcal
L^{\mathrm{n}}(\Theta^+_-)
\cup
\mathcal
L^{\mathrm{n}}(\Theta^-_+)
\cup
\mathcal
L^{\mathrm{n}}(\Theta^-_-)
\cup 
\bigl\{
\lambda^{\mathrm
{rl}}_{\vec 0_1},
\lambda^{\mathrm
{rl}}_{\vec 0_2}
\lambda^{\mathrm
{rl}}_{\overrightarrow \infty_1}
\lambda^{\mathrm
{rl}}_{\overrightarrow \infty_2}
\bigr\},
\end{equation*}
where the oriented slopes
$\vec 0_1$,
$\vec 0_2$,
$\overrightarrow \infty_1$,
$\overrightarrow \infty_2$
correspond to the unoriented
slopes 0 and
$\infty$, is compact. 
\end{proof}

Let
$\mathcal L_0^{\mathrm{ocr}}(T)$
denote the space of all oriented
chain-recurrent geodesic
laminations in the interior of
$T$. There is a natural map 
$\pi \co \mathcal
L_0^{\mathrm{ocr}}(T)
\rightarrow
\mathcal L_0^{\mathrm{cr}}(T)$
defined by forgetting the
orientation. This is a 2--fold
covering map, since
Proposition~\ref
{prop:ChainRecPunctTor} shows that
every chain-recurrent geodesic
lamination is orientable
and connected. 

By Proposition~\ref
{prop:NoSink} and by
inspection in Proposition~\ref
{prop:NonChainRecPunctTor}, each
$\lambda
\in 
\mathcal L_0^{\mathrm{n}}(T)$
uniquely decomposes as the union
of a sink $\sigma(\lambda) \in
\mathcal L_0^{\mathrm{ocr}}(T)$
and of one infinite isolated leaf
whose two ends spiral around
$\sigma(\lambda)$ in the
direction of its orientation.
This defines a map
$\sigma\co \mathcal
L^{\mathrm{n}}(T)
\rightarrow
\mathcal L_0^{\mathrm{ocr}}(T)$.

\begin{lemma}
The map $\sigma\co \mathcal
L^{\mathrm{n}}(T)
\rightarrow
\mathcal L_0^{\mathrm{ocr}}(T)$
is continuous.
\end{lemma}
\begin{proof}
It clearly suffices to show that
the restriction of
$\pi\circ\sigma$ to each subset
$\mathcal
L^{\mathrm{n}}(\Theta^\pm_\pm)$
is continuous. For 
$\lambda \in\mathcal
L^{\mathrm{n}}(\Theta^\pm_\pm)$, we
observed in the proof of
Lemma~\ref{lem:NonRecCompact} that
$\pi\circ\sigma(\lambda) \in
\mathcal L_0^{\mathrm{cr}}(T)$ is
obtained by removing from
$\lambda$ the infinite isolated
leaf that realizes the edge path 
$\left\langle e
\right\rangle$. It follows that
the restriction of
$\pi\circ\sigma$ to
 $\mathcal
L^{\mathrm{n}}(\Theta^\pm_\pm)$
is distance non-increasing for the
metric $d_{\Theta_\pm^\pm}$, and
is therefore continuous. 
\end{proof}

\begin{theorem}
\label{thm:TopNonChainRec}
The space 
$\mathcal L^{\mathrm
n}_0(T)=\mathcal L_0(T)-\mathcal
L_0^{\mathrm {cr}}(T)$
is
homeomorphic to the subspace $K\cup
L_3$ of
$\mathbb R \cup \{\infty\}$ union
of the standard Cantor set
$K\subset [0,1] \subset \mathbb R$
and of a countable set
$L_3$ consisting of exactly  $3$
isolated points in each component of
$\mathbb R \cup
\{\infty\} -K$. 
\end{theorem}

\begin{proof} 
In Theorem~\ref{thm:TopL_0cr}, we
constructed a homeomorphism
$\varphi$ from the space
$\mathcal L_0^{\mathrm{cr}}(T)$
of the chain-recurrent geodesic
laminations to the union
$K\cup L_1$ of $K$ and of
a family $L_1$ of isolated points
consisting of the point
$\infty$ and of exactly
one point in each component of
$[0,1] -K$. Select the
homeomorphism $\varphi$ so that
it satisfies the conditions of
Proposition~\ref{prop:TopL_0Bis}.
In particular, $\varphi$
establishes an order-preserving
one-to-one correspondence between
$\mathbb Q$ and the set of
components of
$[0,1]-K$, by associating to
$s\in \mathbb Q$ the component of
$[0,1]-K$ that contains the
isolated point
$\varphi(\lambda_s)$. 

Consider the 2--fold covering map 
$\mathcal L_0^{\mathrm{ocr}}(T)
\rightarrow
\mathcal L_0^{\mathrm{cr}}(T)$
defined by forgetting the
orientation.  For every
component $I$ of $\mathbb R \cup
\{\infty\} -\varphi \left(
\mathcal L_0^{\mathrm{cr}}(T)
\right)$, it follows from
Proposition~\ref
{prop:TopL_0Bis}(iii) that an
orientation of the geodesic
orientation corresponding under
$\varphi$ to one end point of $I$
uniquely determines an
orientation of the geodesic
lamination corresponding to the
other end point. We can
consequently lift $\varphi\co
\mathcal L_0^{\mathrm{cr}}(T)
\rightarrow \mathbb R \cup
\{\infty\}$ to a
continuous map $\widetilde
\varphi \co
\mathcal L_0^{\mathrm{ocr}}(T)
\rightarrow
\mathbb S^1$, where we denote by
$\mathbb S^1$ the circle that is
the 2--fold covering of $\mathbb R
\cup \{\infty\}$. (Of course,
$\mathbb S^1$ is homeomorphic to
$\mathbb R\cup \{\infty\}$, but
we prefer to use a  different
letter to emphasize the
distinction). 

Let $\widetilde K$ denote the
preimage of $K$ in $\mathbb S^1$.
Pick any family $L_3$ of isolated
points in $\mathbb S^1 - \widetilde
K$ such that, for every
component $I$ of $\mathbb S^1-
\widetilde K$, the intersection
$I\cap L_3$ consists of exactly
3 points $x_I^{\mathrm r}$,
$x_I^{\mathrm l}$ and
$x_I^{\mathrm {rl}}$. Note that
each component $I$ of $\mathbb
S^1 -\widetilde K$ is now indexed
by an oriented slope
$\vec s$, namely $I$ is
the component $I_{\vec
s}$
 containing the image under
$\widetilde \varphi$ of
the oriented simple closed
geodesic
$\lambda_{\vec s}$ of
oriented slope $\vec
s$, and whose boundary points are
consequently $\widetilde \varphi
(\lambda_{\vec
s}^-)$ and $\widetilde \varphi
(\lambda_{\vec
s}^+)$ by
Proposition~\ref{prop:TopL_0Bis}.

We now define a map $\psi\co 
\mathcal L_0^{\mathrm{n}}(T)
\rightarrow \widetilde K \cup L_3$ as
follows. If the sink
$\sigma(\lambda)
\in
\mathcal L_0^{\mathrm{ocr}}(T)$
of $\lambda
\in 
\mathcal L_0^{\mathrm{n}}(T)$
 is not a simple
closed geodesic, define
$\psi(\lambda)$ as $\widetilde
\varphi(\sigma(\lambda)) \in
\widetilde K $. Otherwise,
$\sigma(\lambda)$ is the oriented
simple closed geodesic
$\lambda_{\vec s}$ for
some oriented slope
$\vec s$, and
$\lambda$ is the geodesic
lamination 
$\lambda_{\vec
s}^{\mathrm{r}}$, 
$\lambda_{\vec
s}^{\mathrm{l}}$ or
$\lambda_{\vec
s}^{\mathrm{rl}}$ with the
notation of Proposition~\ref
{prop:NonChainRecPunctTor};
in this case, define $\psi
(\lambda)$ as the isolated point
$x^{\mathrm{r}}
_{I_{\vec  s}}$, 
$x^{\mathrm{l}}
_{I_{\vec  s}}$ or
$x^{\mathrm{rl}}
_{I_{\vec  s}} \in
L_3$, respectively. 

We will show that $\psi\co 
\mathcal L_0^{\mathrm{n}}(T)
\rightarrow
\widetilde K \cup L_3$ is a
homeomorphism.

Because $\varphi$ is a
homeomorphism, $\widetilde
\varphi$ is injective and it
immediately follows from the
construction that $\psi$ is a
bijection.

To prove that 
$\psi\co 
\mathcal L_0^{\mathrm{n}}(T)
\rightarrow
\widetilde K \cup L_3$
is continuous, we will show that
for every sequence
$\alpha_n
\in
\mathcal L_0^{\mathrm{n}}(T)$, $n
\in \mathbb N$,
converging to $\lambda \in 
\mathcal L_0^{\mathrm{n}}(T)$, the
sequence $\psi(\alpha_n)$ admits
a subsequence converging to
$\psi(\lambda)$. Passing to a
subsequence if necessary, we can
assume that, either the
sink $\sigma(\alpha_n)
\in \mathcal
L_0^{\mathrm{ocr}}(T)$ is a
closed geodesic for every $n$, or
it is a closed geodesic for no
$n$.
If the
$\sigma(\alpha_n)$ are not closed
geodesics, then
$\psi(\alpha_n)=\widetilde
\varphi (\sigma(\alpha_n))$
converges to $ \widetilde
\varphi(\sigma(\lambda)=
\psi(\lambda)$ by continuity of
$\widetilde \varphi$ and
$\sigma$, and we are done. We can
consequently assume that each
$\sigma(\alpha_n) $ is an
oriented closed geodesic
$\lambda_{\vec s_n}$
of oriented rational slope
$\vec s_n$. 

If the
limit $\sigma(\lambda)$ of
$\sigma(\alpha_n)$ is a closed
geodesic $\lambda_{\vec
{s}}$, then it is isolated in 
$\mathcal
L_0^{\mathrm{ocr}}(T)$ and
$\sigma(\alpha_n) =
\lambda_{\vec {s}}$
for $n$ large enough. It follows
that 
$\alpha_n=
\lambda_{\vec
s}^{\mathrm{r}}$, 
$\lambda_{\vec
s}^{\mathrm{l}}$ or
$\lambda_{\vec
s}^{\mathrm{rl}}$, and therefore
that the converging sequence
$\alpha_n$ is eventually
constant, equal to its limit
$\lambda$. In particular,
$\psi(\alpha_n)$ converges to
$\psi(\lambda)$. 

If $\sigma(\lambda)$ is not a
closed geodesic, the sequence
$\vec s_n$ has
no constant subsequence. It
follows that the length of the
component of
$\mathbb S^1 - \widetilde K$
containing $\psi(\alpha_n)$ tends
to 0 as $n$ tends to $\infty$.
Since this component also
contains the point $\widetilde
\varphi\circ\sigma (\alpha_n)$ by
construction of $\psi$, we
conclude that the sequence
$\psi(\alpha_n)$ converges
to the limit of 
$\widetilde
\varphi\circ\sigma (\alpha_n)$,
namely to 
$\widetilde
\varphi\circ\sigma (\lambda) =
\psi(\lambda)$ by continuity of
$\widetilde\varphi$ and $\sigma$.

This concludes the proof that the
bijection $\psi\co 
\mathcal L_0^{\mathrm{n}}(T)
\rightarrow
\widetilde K \cup L_3$ is continuous.
Because $\mathcal
L_0^{\mathrm{n}}(T)$ is compact
by Lemma~\ref{lem:NonRecCompact},
it follows that
$\psi$ is a homeomorphism. Since 
there is a homeomorphism
from $\mathbb S^1$ to $\mathbb R
\cup \{\infty\}$ sending
$\widetilde K$ to $K$, this 
concludes the proof of
Theorem~\ref
{thm:TopNonChainRec}. 
\end{proof}

\section{The Hausdorff dimension of
$\mathcal L_0(T)$ for the
once-punctured torus}

\begin{theorem}
\label{thm:DimL_0}
The space $\left( \mathcal
L_0(T) , d_{\log}
\right)$ has Hausdorff dimension
$2$, and its $2$--dimensional
Hausdorff measure is equal to
$0$. 
\end{theorem}
\begin{proof}
Since $\mathcal L_0(T)$ contains
$\mathcal L_0^{\mathrm {cr}}(T)$,
which has Hausdorff dimension 2 by
Theorem~\ref{thm:DimL_0cr}, the
Hausdorff dimension of $\left(
\mathcal L_0(T) , d_{\log}
\right)$  is at least 2.
Therefore, it suffices to show
that its 2--dimensional Hausdorff
measure is equal to 0. 

Fix $\varepsilon>0$ and $r>0$. By
Theorem~\ref{thm:DimL_0cr},
$\left( \mathcal
L_0^{\mathrm{cr}}(T), d_{\log}
\right)$ has 2--dimensional
measure 0. In particular,
the subset $\mathcal
L_0^{\mathrm{r}}(T)$ consisting
of all recurrent geodesic
laminations also has
2--dimensional Hausdorff measure
0. Therefore, we can cover
$\mathcal L_0^{\mathrm{r}}(T)$
by a family of $d_{\log}$--balls
$B(\lambda_i, r_i)$, $i\in I$,
with
$\lambda_i \in \mathcal
L_0^{\mathrm r}(T)$, $r_i<r$ and
$\sum_{i \in I} r_i^2<\varepsilon$. 

Let $r_0$ and $c_0$ be
the constants of
Lemma~\ref{lem:ConvNonRecurPart},
and assume $r
<r_0$ without loss of
generality. For each of the above
balls $B(\lambda_i, r_i)$,
consider the balls $B(\lambda',
r')$ where $\lambda'$
contains $\lambda_i$ and where
$r'=c_0r_i$. By
Propositions~\ref
{prop:ChainRecPunctTor} and
\ref{prop:NonChainRecPunctTor},
there are at most 11 such
$\lambda'$. By Lemma~\ref
{lem:ConvNonRecurPart}, we can
therefore cover the whole space
$\mathcal L_0(T)$ by a family of
balls $B(\lambda'_j,
r'_j)$, $j\in J$, such that
$r'_j<c_0 r$ and
$\sum_{j\in J} (r'_j)^2 <11 c_0^2
\varepsilon$. 
Since this holds for every
$r<r_0$ and every
$\varepsilon$, this proves that 
 $\left(
\mathcal L_0(T) , d_{\log}
\right)$  has 2--dimensional
Hausdorff measure  0. 
\end{proof}

\section{The 4--times-punctured
sphere}
\label{sect:PunctSphere}

We now consider the case where
the surface $S$ is the 4--times-punctured sphere. The analysis is
very similar to that of the once-punctured torus, and we will only
sketch the arguments. 

Consider the group $\Gamma$ of
diffeomorphisms of $\mathbb R^2 -
\mathbb Z^2$  consisting of all
rotations of $\pi$ around the
points of the lattice $\mathbb
Z^2$, and of all translations by
the elements of the sublattice
$\left(2\mathbb Z\right)^2$. The
quotient space $\left(\mathbb R^2
-
\mathbb Z^2\right)/\Gamma$ is
diffeomorphic to the interior of
the 4--times-punctured sphere
$S$.

What makes the 4--times-punctured sphere so similar to
the once-punctured torus is
that, in both cases, simple
closed geodesics in the
interior of the surface are
characterized by their slope. As
in the case of the once-punctured
torus, every straight line with
rational slope
$\frac pq \in \mathbb Q\cup
\{\infty\}$ projects to a simple
closed curve in $S$. Conversely,
every simple closed curve that is
not isotopic to a boundary
component is obtained in this way.
This establishes a one-to-one
correspondence between the set
$\mathcal S(S)$ of simple closed
geodesics in the interior of the
4--times-punctured sphere $S$ and
the set of rational slopes 
$\frac pq \in \mathbb Q\cup
\{\infty\}$. 

\begin{figure}[ht]
\centerline{
\includegraphics
{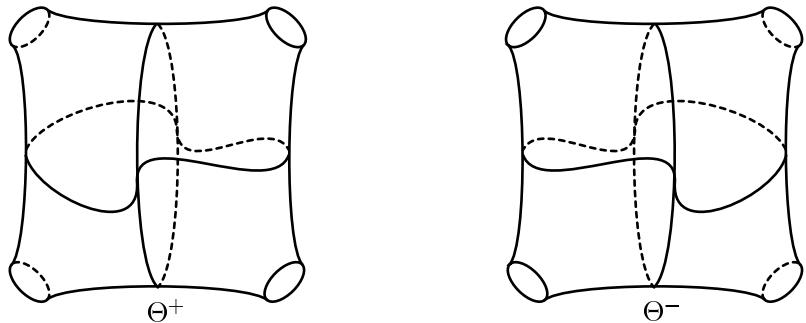}}
\caption{The train tracks
$\Theta^+$ and $\Theta^-$ on the
4--times-punctured sphere $S$}
\label
{pict:TrainTracks4PunctSphere}
\end{figure}

Consider the train tracks
$\Theta^+$ and $\Theta^-$ shown
on Figure~\ref
{pict:TrainTracks4PunctSphere}.
For the appropriate identification
between $\mathbb R^2/\Gamma$ and 
the interior of
$S$, the preimage of
$\Theta^+$ is exactly the train
track which
already appeared in
Figure~\ref{pict:R2-Z2} for the
once-punctured torus. In
particular, every simple closed
geodesic $\lambda \in \mathcal
S(S)$ with non-negative slope
$\frac pq \in [0,\infty]\cap
\mathbb Q$ is weakly carried by
$\Theta^+$, and every simple
closed geodesic with non-positive
slope is weakly carried by
$\Theta^-$. 

The key estimate is the
following analog of
Proposition~\ref
{prop:d_phiAndSlopes}, whose
proof very closely follows that
of that first result.

\begin{proposition}
\label
{prop:d_phiAndSlopes4PuncSphere}
On the $4$--times-punctured
sphere $S$, let the simple
closed geodesics
$\lambda$, $\lambda'\in\mathcal
S(S)$ have slopes
$\frac pq$, $\frac {p'}{q'} \in
\mathbb Q\cup \{\infty\}$ with
$0\leq \frac pq  <\frac
{p'}{q'}\leq
\infty$. Then 
\begin{equation*} 
\textstyle
d_{\Theta^+}
(\lambda, \lambda') = 
\max \Bigl\{ \frac1{\lvert
p''\rvert + \lvert q''\rvert};
 \frac pq \leq\frac {p''}{q''}
\leq\frac {p'}{q'}
\Bigr\}.
\end{equation*} 
\end{proposition}

As in the case of the once-punctured torus, the combination
of Proposition~\ref
{prop:d_phiAndSlopes4PuncSphere}
and of Proposition~\ref
{prop:TopQhat} in the
Appendix provide a homeomorphism
between the space $\mathcal
L_0^{\mathrm{cr}}(S)$ of
chain-recurrent geodesic laminations
and the subspace $K\cup L_1$ of $
\mathbb R
\cup
\{\infty\}$  obtained by adding
to the standard middle third
Cantor set
$K\subset [0,1]$ a family $L_1$ of
isolated points, consisting of
the point
$\infty$ and of exactly one
isolated point
 in each component of
$[0,1] - K$.

The space 
$\mathcal
L_0^{\mathrm{ncr}}(S)$
of non-chain-recurrent
geodesic laminations is much
simpler for the 4--times-punctured sphere $S$ than for the
once-punctured torus $T$. To see
this, we need to analyze the
topology of geodesic laminations
in $S$.

We begin by borrowing from
\cite{PenHar} or \cite{FatLauPoe}
the classification of measured
geodesic laminations. For the
identifications $T\cong
\bigl( \mathbb R^2 -
\mathbb Z^2 \bigr)/\mathbb Z^2$
and
$S\cong
\bigl( \mathbb R^2 -
\mathbb Z^2 \bigr)/\Gamma$, it can be
shown that, for every slope $s\in
\mathbb R \cup\{\infty\}$, there
is a geodesic lamination $\mu_s$
in the interior of the 4--times-punctured sphere $S$ whose
preimage to $ \mathbb R^2 -
\mathbb Z^2 $ coincides with the
preimage of the geodesic
lamination $\lambda_s\in \mathcal
L(T)$ of Section~\ref
{sect:TopGeodLamPuncTor}.

\begin{proposition}
\label{prop:Rec4PuncSph}
Every recurrent geodesic
lamination in the interior of the
$4$--times-punctured sphere is of
the form
$\mu_s$ for some $s\in \mathbb
R\cup \{\infty\}$. When $s$ is
rational, $\mu_s$ is a simple
closed geodesic, and the
completion of its complement
consists of two twice-punctured
disks. Otherwise, $\mu_s$ has
uncountably many leaves and the
completion of its complement
consists of four once-punctured
monogons, each with one spike. 
\end{proposition}

In particular, the completion of
the complement of each $\mu_s$
again contains only finitely many
simple geodesics. A case-by-case
analysis then provides the
following two statements.

\begin{proposition}
\label{prop:ChainRec4PuncSph}
The chain-recurrent geodesic
laminations in the interior of
the $4$--times-punctured sphere
$S$ fall into the following
categories:
\begin{enumerate}
\item The recurrent geodesic
lamination $\mu_s$, with $s\in
\mathbb R\cup \{\infty\}$.

\item The union $\mu_s^+$ of the
simple closed geodesic $\mu_s$,
with rational slope $s \in
\mathbb Q\cup \{\infty\}$, and of
two infinite isolated leaves, one
in each component of $S-\mu_s$;
for an arbitrary orientation of
$\mu_s$, the two ends of the leaf
in the left component of $S-\mu_s$
spiral along $\mu_s$ in the
direction of the orientation, and
the ends of the leaf in the right
component of $S-\mu_s$ spiral
along $\mu_s$ in the opposite
direction.

\item The union $\mu_s^-$ of the
simple closed geodesic $\mu_s$,
with rational slope $s \in
\mathbb Q\cup \{\infty\}$, and of
two infinite isolated leaves, one
in each component of $S-\mu_s$;
for an arbitrary orientation of
$\mu_s$, the two ends of the leaf
in the right component of
$S-\mu_s$ spirals along $\mu_s$
in the direction of the
orientation, and the ends of the
leaf in the left component of
$S-\mu_s$ spiral along $\mu_s$ in
the opposite direction.
\end{enumerate}
\end{proposition}

Note that the $\mu_s$ with
irrational slopes $s$, as well
as the $\mu_s^+$ and $\mu_s^-$
(with rational slopes) are
non-orientable. In particular, by
Proposition~\ref
{prop:NoSink}, every
non-chain-recurrent geodesic
lamination is obtained by adding to
a simple closed geodesic $\mu_s$,
with rational slope $s\in \mathbb Q
\cup\{\infty\}$, a certain number
of infinite isolated leaves whose
ends all spiral along $\mu_s$ in
the same direction. Looking at
possibilities, one concludes:

\begin{proposition}
\label{prop:NonChainRec4PuncSph}
Any non-chain-recurrent geodesic
lamination in the interior of
the $4$--times-punctured sphere
is the union of a
closed geodesic $\mu_s$, with
rational slope $s\in \mathbb Q
\cup\{\infty\}$, and of $1$ or
$2$ infinite isolated leaves
whose ends all spiral along
$\mu_s$ in the same direction. 
A given simple closed geodesic
$\mu_s$ is contained in exactly
$6$ such non-chain-recurrent
geodesic laminations.
\end{proposition}

As in the proof of Theorem~\ref
{thm:TopNonChainRec}, we can then
use Propositions~\ref
{prop:ChainRec4PuncSph} and \ref
{prop:NonChainRec4PuncSph}
to push our analysis of the
topology of $\mathcal
L_0^{\mathrm{cr}}(S)$ to $\mathcal
L_0(S)$. 

\begin{theorem}
For the $4$--times-punctured
sphere $S$, the space $\mathcal
L_0(S)$ is homeomorphic to the
subspace $K \cup L_7$ of
$\mathbb R
\cup
\{\infty\}$ union of the standard
middle third Cantor set $K
\subset [0,1]$ and of a set $L_7$ of
isolated points consisting of
exactly $7$ points in each
component of $\mathbb R \cup
\{\infty\} -K$. The homeomorphism
can be chosen so that the set 
$\mathcal S(S) \subset \mathcal
L_0(S)$ of all simple
closed geodesics corresponds to a
subset $L_1\subset L_7$
consisting of exactly $1$ point
in each component of $\mathbb R
\cup
\{\infty\} -K$. The closure of
$\mathcal S(S)$, namely the
space $\mathcal
L_0^{\mathrm{cr}}(S)$ of
chain-recurrent geodesic
laminations, then corresponds to
the union $K
\cup L_1$. Its
complement, the space  $\mathcal
L_0^{\mathrm{n}}(S)$ of
non-chain-recurrent geodesic
laminations, is countable.
\end{theorem}

More precisely, each component $I$
of 
$\mathbb R
\cup
\{\infty\} -K$ is indexed by the
rational slope $s \in \mathbb Q
\cup \{\infty\}$ of the simple
closed geodesic $\mu_s$ whose
image is contained in $I$. The
end points of the interval $I$
then correspond to the
chain-recurrent geodesic laminations
$\mu_s^\pm$, and the
intersection of the interior of
$I$ with $\widehat X$
corresponds to $\mu_s$ and to
the 6 non-chain-recurrent
geodesic laminations containing
it.

The Hausdorff dimension and
measure of $\left (\mathcal
L_0(S), d_{\log} \right)$ are
obtained by combining
Propositions~\ref
{prop:d_thetaEqd_log}, 
\ref
{prop:d_phiAndSlopes4PuncSphere},
\ref
{prop:dim2} and the fact that
$\mathcal L_0^{\mathrm{ncr}}(S)$ 
is countable. 

\begin{theorem} 
\label{thm:HausDim4PunSph}
For the
$4$--times-punctured sphere $S$,
the Hausdorff dimension of the
metric space $\left (\mathcal
L_0(S), d_{\log} \right)$ is
equal to
$2$. Its $2$--dimensional
Hausdorff measure is equal to
$0$.
\end{theorem}

\section{Very small surfaces}
\label{sect:VerySmall}

Having considered the once-punctured torus
or the 4--times-punctured
sphere, we may wonder
about surfaces of lower
complexity. 
 Geodesic
laminations on a surface $S$ make
sense only when
$S$ admits a metric of negative
curvature for which the boundary is
totally geodesic, namely when the
Euler characteristic $\chi(S)$ is
negative. This leaves the 3--times-punctured
sphere, the twice-punctured projective plane and the
once-punctured Klein bottle. 

We
will see that the geodesic lamination
spaces of these surfaces are
relatively trivial, thereby
justifying the emphasis of this
paper on the once-punctured torus
and on the 4--times-punctured
sphere.

\begin{proposition} 
\label{prop:FiniteSpaces} If the
surface $S$ is the 3--times-punctured
sphere or the twice-punctured
projective plane,
then the space $\mathcal L(S)$ of
geodesic laminations on $S$ is
finite.
\end{proposition}

\begin{proof} Each of these
surfaces has only finitely many
homotopy classes of simple closed
curves. It follows that they
contain only finitely many
multicurves and therefore that
every chain-recurrent geodesic
lamination is a multicurve. In
particular, every recurrent
geodesic lamination is a
multicurve, since
it is chain-recurrent by
Proposition~\ref{prop:NoSink}.

If $S$ is the 3--times-punctured
sphere, each multicurve is in
addition contained in the boundary
$\partial S$. There  are only
finitely many simple arcs 
$a\subset S$ with $\partial
a\subset \partial S$, modulo
homotopy keeping $\partial a$ in
$\partial S$. It easily follows
that there are only finitely many 
infinite simple geodesics in $S$,
spiralling along boundary
components. This implies that the
3--times-punctured sphere $S$
contains only finitely many geodesic
lamination. 

When $S$ is the twice-punctured
projective plane, splitting $S$ open
along a multicurve $\lambda_1$
produces a twice-punctured
projective plane (when
$\lambda_1\subset \partial S$) or a
3--times-punctured sphere. Again,
the twice-punctured projective plane
contains only finitely many homotopy
classes of simple arcs, relative to
the boundary. In both cases, it
follows that
$\lambda_1$ can be extended to a
finite number of geodesic
laminations. Since $S$ contains a
finite number of multicurves,
$\mathcal L(S)$ is finite for the
twice-punctured projective plane
$S$. 
\end{proof}

For the once-punctured Klein bottle,
we restrict attention to the closed
subspace $\mathcal L_0(S)\subset
\mathcal L(S)$ consisting of those
geodesic laminations that are
contained in the interior of $S$. As
indicated in the introduction, this
is essentially an exposition choice
as the results easily extend to the
whose space
$\mathcal L(S)$, but at the expense
of more cases to consider.

\begin{proposition} 
\label{prop:KleinBottle} If $S$ is
the once-punctured Klein Bottle,
then the space $\mathcal L_0(S)$ of
geodesic laminations  in the
interior of $S$ is countable
infinite\footnote{And not finite, as
erroneously stated in \cite{ZhuBon1}
and in an earlier version of the
current paper!}. All of its points
are isolated, with the exception of
six limit points. The closure
$\mathcal L_0^{\mathrm{cr}}(S)$ of
the set of multicurves consists of
infinitely many isolated points and
of two limit points; these two limit
points of $\mathcal
L_0^{\mathrm{cr}}(S)$ are also limit
points of $\mathcal
L_0(S)-\mathcal
L_0^{\mathrm{cr}}(S)$.
\end{proposition}

\begin{proof} Up to homotopy, the
interior of the once-punctured Klein
bottle $S$ contains only one
orientation-preserving simple closed
curve that is essential, in the
sense that it is not parallel to the
boundary and that it bounds neither a
disk nor a M\"obius strip. It follows
that the interior of $S$ contains
only one orientation-preserving
simple closed geodesic
$\lambda_\infty$ in the interior of
$S$. 

There are infinitely many
homotopy classes
of simple, orientation-reversing,
closed curves, but these are easily
classified. The corresponding
orientation-reversing simple closed
geodesics can be listed as
$\lambda_n$, $n\in \mathbb Z$, is
such a way that each $\lambda_n$ is
homotopic to $T^n(\lambda_0)$ where
$T$ denotes the Dehn twist around
$\lambda_\infty$ (well-defined once
we fix an orientation on a
neighborhood of the
orientation-preserving simple closed
geodesic $\lambda_\infty$, and
once we arbitrarily decide which
simple closed geodesic will be
called $\lambda_0$). In particular,
each $\lambda_n$ meets
$\lambda_\infty$ in exactly one
point. 

In addition, $\lambda_m$ is disjoint
from $\lambda_n$ exactly when
$m=n\pm1$.

As $n$ tends to $+\infty$, the
simple closed geodesic $\lambda_n$
converges to a geodesic lamination
$\lambda_\infty^-$, union of
$\lambda_\infty$ and of an infinite
simple geodesic whose ends spiral on
each side of $\lambda_\infty$, in
opposite direction. As $n$ tends to
$-\infty$, $\lambda_n$ converges to
a geodesic lamination
$\lambda_\infty^+$ of the same type,
but with opposite spiralling
directions on each side of
$\lambda_\infty$. 

It follows that chain-recurrent
geodesic laminations are of the
following four possible types:
$\lambda_\infty$, $\lambda_n$,
$\lambda_n \cup \lambda_{n+1}$,
$\lambda_\infty^+$ or
$\lambda_\infty^-$. In particular,
the subspace $\mathcal L_0^{\mathrm
{cr}}(S)$ consisting of all
chain-recurrent geodesic lamination
is countable infinite. All of its
points are isolated, with the
exception of two limit points
corresponding to
$\lambda_\infty^\pm$. 

From the above list, we see that
splitting $S$ along a recurrent
geodesic lamination $\lambda$
produces a 3--times-punctured
sphere, or a once-punctured bigon,  
or a twice-punctured
projective plane. By
Proposition~\ref{prop:FiniteSpaces},
it follows that such a recurrent
geodesic lamination can be enlarged
to finitely many geodesic
laminations. By
Proposition~\ref{thm:DynGeodLam},
this implies that the space
$\mathcal L_0(S)$ is countable.

A limit point
of  $\mathcal
L_0(S)$ must contain a limit point
of $\mathcal
L_0^{\mathrm{cr}}(S)$, by
consideration of the recurrent parts
of geodesic laminations converging
to that limit. Therefore, a limit point
of  $\mathcal
L_0(S)$ must contain
$\lambda_\infty^+$ or
$\lambda_\infty^-$. The complement
of the geodesic lamination
$\lambda_\infty^\pm$ is a punctured
bigon. As such, this complement
contains exactly two simple
geodesics, each of which has its two
ends converging towards one of the
spikes of the bigon. It follows that
$\lambda_\infty^\pm$ can be enlarged
to exactly two geodesic laminations.
By considering suitable
enlargements of $\lambda_n$ and
letting $n$ tend to $\pm\infty$, one
easily sees that these enlarged
geodesic laminations are indeed limit
points of $\mathcal L_0(S)$.
Therefore, $\mathcal L_0(S)$ has
exactly 6 limit points.    

Similarly, the two elements 
$\lambda_\infty^\pm$ are  limit
points of $\mathcal
L_0(S)-\mathcal
L_0^{\mathrm{cr}}(S)$.
\end{proof}

Since the geodesic lamination spaces
$\mathcal L_0(S)$ of the 3--times-punctured sphere, the twice-punctured projective plane and the
once-punctured Klein bottle are
countable, their Hausdorff dimension
is of course 0 for any metric, and
in particular for the metric
$d_{\log}$.

\section*{Appendix}
\addcontentsline{toc}{section}{Appendix}

In this appendix, we study the
space
$\mathbb Q \cup \{\infty\}$ with
the (ultra)metric $d$ defined by 
\begin{equation*} 
\textstyle 
d\bigl(
\frac pq,
\frac{p'}{q'}\bigr) =\max
\Bigl\{  \frac 1{\left\vert
p'' \right\vert +q''};
\frac pq
\leq
\frac{p''}{q''} \leq
\frac{p'}{q'} \Bigr \}
\end{equation*}
 when $\frac pq <\frac {p'}{q'}$.
More precisely, we will study the
completion
$\widehat {\mathbb Q}$ of $\mathbb
Q \cup \{\infty\}$ for this
metric.

By convention, whenever we
consider a rational number $\frac
pq \in
\mathbb Q \cup \{\infty\}$, we
always implicitly assume that $p$
and $q$ are integer and coprime,
and that $q\geq 0$. In
particular, $0=\frac 01$ and
$\infty = \frac 10=\frac {-1}0$.

Recall that the
completion $\widehat {\mathbb Q}$
can be defined as the set of
equivalence classes of Cauchy
sequences in
$\mathbb Q \cup \{\infty\}$,
where two Cauchy sequences are
equivalent when their union is a
Cauchy sequence. We consequently
need to analyze Cauchy sequences
in
$(\mathbb Q \cup \{\infty\},d)$.

\begin{lemma}
\label{lem:Cauchy}
A sequence
$\bigl(\frac{p_n}{q_n}\bigr)_{n
\in
\mathbb N}$ in $\left( \mathbb Q
\cup
\{\infty\}, d\right)$ is Cauchy
if and only if one of the
following holds:
\begin{enumerate}
\item 
\label{item1:Cauchy}
The sequence
converges to an irrational number
for the usual topology of
$\mathbb R \cup
\{\infty\}$.
\item 
\label{item2:Cauchy}
The sequence converges to a
rational number $\frac pq \in
\mathbb Q \cup \{\infty\}$ for the
usual topology of $\mathbb R \cup
\{\infty\}$ and, for $n$
sufficiently large,
$\frac{p_n}{q_n}$ stays on one
side of $\frac pq$.
\item 
\label{item3:Cauchy}
For $n$ sufficiently large,
$\frac{p_n}{q_n}$ is equal to a
fixed rational number $\frac pq
\in \mathbb Q \cup \{\infty\}$. 
\end{enumerate}
\end{lemma}

\begin{proof} First assume that
the sequence 
$\bigl(\frac{p_n}{q_n}\bigr)_{n
\in
\mathbb N}$ is Cauchy for the
metric $d$.

By compactness, $\bigl(\frac{p_n}{q_n}\bigr)_{n
\in
\mathbb N}$ admits a subsequence
which converges for the usual
topology of $\mathbb R \cup
\{\infty\}$. If two
subsequences of $\bigl(\frac{p_n}{q_n}\bigr)_{n
\in
\mathbb N}$ had different limits $l\not=l'$ then, for an
arbitrary rational number $\frac pq$  between $l$ and $l'$, we
would have
$d\bigl(\frac {p_m}{q_m}, \frac
{p_n}{q_n}\bigr)\geq
\frac1{\left|p\right|+q}$ whenever $\frac {p_m}{q_m}$ is close
to $l$ and $\frac {p_n}{q_n}$ is close to $l'$, contradicting
the fact that the sequence is Cauchy. We conclude
that $\bigl(\frac{p_n}{q_n}\bigr)_{n
\in
\mathbb N}$  must have a limit
$x$ for the usual topology of
$\mathbb R \cup
\{\infty\}$. 

If the limit $x$ is a rational
number $\frac pq$, then
$d\bigl(\frac {p_m}{q_m}, \frac
{p_n}{q_n}\bigr)\geq
\frac1{\left|p\right|+q}$ if
$\frac {p_m}{q_m}\leq \frac pq
\leq \frac {p_n}{q_n}$, unless 
$\frac {p_m}{q_m}= \frac
pq= \frac {p_n}{q_n}$. This
proves that, either 
$\frac {p_m}{q_m}$
stays on one side of $\frac pq$
for $n$ sufficiently large, or 
$\frac {p_m}{q_m}= \frac pq$
for $n$ sufficiently large. 

Consequently, if the sequence 
$\bigl(\frac{p_n}{q_n}\bigr)_{n
\in
\mathbb N}$ is Cauchy for $d$,
then it satisfies Conditions~%
\ref{item1:Cauchy}--%
\ref{item3:Cauchy} of the lemma.

Conversely, if $\bigl(\frac{p_n}{q_n}\bigr)_{n
\in \mathbb N}$ satisfies
Conditions~\eqref{item1:Cauchy}--\eqref{item3:Cauchy}, one easily
sees that it is Cauchy for $d$.
\end{proof}

\begin{corollary} 
\label{cor:IrrationalPoints}The
completion
$\widehat {\mathbb Q}$ contains a
subset $A$ which is isometric to
$\mathbb R -\mathbb Q$ endowed
with the metric $d$ defined by 
\begin{equation*}
d(x,y) = \max 
\Bigl\{
\textstyle
\frac1{\left|p\right|+q}; x<\frac
pq <y \Bigr\}
\end{equation*}
for every $x<y$ in $\mathbb R
-\mathbb Q$. In addition, the
complement $\widehat {\mathbb
Q} - A$ is countable. 
\end{corollary}
\begin{proof} The fact
that every Cauchy sequence
in $\left( \mathbb Q
\cup
\{\infty\}, d\right)$
admits a limit for the
usual topology of $\mathbb
R \cup \{\infty\}$ defines
a continuous map $\pi\co
\widehat {\mathbb Q}
\rightarrow \mathbb R \cup
\{\infty\}$. By
Lemma~\ref{lem:Cauchy}, the
preimage of an irrational point under $\pi$ consists of a
single point of $\widehat
{\mathbb Q}$. It follows that the
restriction of $\pi$ induces a
bijection between $A = \pi^{-1}
\left( \mathbb R-\mathbb Q\right)$
and
$\mathbb R-\mathbb Q$. In
addition, this restriction of
$\pi$ to $A$ clearly sends the
metric of $\widehat {\mathbb Q}$
to the metric on $\mathbb
R-\mathbb Q$ indicated.

By Lemma~\ref{lem:Cauchy}, the
preimage of a rational point
$\frac pq \in \mathbb Q\cup
\{\infty\}$ under
$\pi$ consists of 3 points in
$\widehat {\mathbb Q}$: one
corresponding to the constant
sequence $\frac pq$, and the
other two corresponding to each
side of
$\frac pq$. It follows that the
complement $\pi^{-1} \left (
\mathbb Q\cup
\{\infty\} \right)$ of $A$ is
countable. 
\end{proof}

We now describe a topological
model for
$\widehat {\mathbb Q}$. Let
$K\subset [0,1] \subset \mathbb
R$ be the standard
middle third Cantor set, and let $X$
be a family of isolated points
consisting of the point
$\infty$ and of exactly one point
 in each component of
$[0,1] - K$.
Note that
$\widehat X=K\cup X$ is equal to the
closure of $X$ in $\mathbb R\cup \{\infty\}$.

\begin{proposition} 
\label{prop:TopQhat}
The
completion $\widehat {\mathbb Q}$ 
of $\mathbb
Q \cup \{\infty\}$ for the
metric $d$ is homeomorphic to the
subspace $\widehat X\subset
\mathbb R\cup \{\infty\}$
described above. 
\end{proposition}

\begin{proof} We first construct
an order-preserving map
$\varphi\co X\rightarrow
\mathbb Q \cup \{\infty\}$,
using the Farey combinatorics of
rational numbers; see for instance
\cite[Section~3.1]{HarWri}. (There are
other possibilities to construct
$\varphi$, but this one seems
prettier). 

In the standard construction of
the Cantor set $K$, the set
$\mathcal I$ of components of
$[0,1]-K$ is
written as an increasing union 
$\mathcal I_n$ of $2^n-1$
disjoint intervals where:
$\mathcal I_1$
consists of the interval $\left]
\frac13,
\frac23\right[$; the set $\mathcal
I_{n+1}$ is obtained from
$\mathcal I_n$ by inserting one
interval of length $\frac1{3^n}$
between any pair of consecutive
intervals of
$\mathcal I_n$, as well as before
the first interval and after the
last one. Since the set
$X-\{\infty\}$ was defined by
picking one point in each
interval of
$\mathcal I$, this provides a
description of
$X$ as an increasing union of
finite sets $X_n$, each with
$2^n$ elements, such that
$X_0=\{\infty\}$ and such that, 
 for $n\geq 1$,
$X_n$ consists of $\infty$
and of exactly one point in each
interval of
$\mathcal I_n$. In
particular, 
$X_{n+1}$ is obtained from $X_n$
by adding one point between each
pair of consecutive elements of
$X_n$. 

The set
$Y=\mathbb Q \cup \{\infty\}$ can
similarly be written as an
increasing union of finite sets
with $2^n$ elements such that:
$Y_0=\bigl\{ \infty \bigr\}$; 
$Y_1=\bigl\{0, \infty
\bigr\}=\bigl\{ \frac01, \frac10
\bigr\}$, $Y_2=\bigl\{\frac{-1}2,
\frac01, \frac12, \frac 10
\bigr\}$; more generally, 
the set
$Y_{n+1}$ is obtained from $Y_n$
by adding the point
$\frac{p+p'}{q+q'}$ between any
two consecutive elements $\frac
pq$, $\frac{p'}{q'}$ of $Y_n$
(counting $\infty$ as both
$\frac10$ and $\frac{-1}0$). See
\cite[Section~3]{HarWri} for a proof
that $Y_n$ contains all the
$\frac pq \in \mathbb Q \cup
\{\infty\}$ with
$\left \vert p\right\vert +q\leq
n+1$, which implies that the
union of the
$Y_n$ is really equal to
$Y=\mathbb Q \cup \{\infty\}$. 

Define $\varphi\co X\rightarrow 
\mathbb Q \cup \{\infty\}$ as
the unique order-preserving map
which sends each $2^n$--element
set $X_n$ to the $2^n$--element
set $Y_n$.

From Lemma~\ref{lem:Cauchy}, one
easily sees that $\varphi$
establishes a one-to-one
correspondence between Cauchy
sequences in $\left( \mathbb Q
\cup
\{\infty\}, d\right)$ and Cauchy
sequences in $X$ (namely
sequences in $X$ which have a
limit in $\widehat X$).
Therefore, $\varphi$ induces a
homeomorphism between the
completion $\widehat {\mathbb
Q}$ of $\left( \mathbb Q
\cup
\{\infty\}, d\right)$ and the
completion $\widehat X$ of $X$. 
\end{proof}

\begin{proposition}  
\label{prop:dim2}
The
completion $\widehat {\mathbb Q}$ 
of $\mathbb
Q \cup \{\infty\}$ for the
metric $d$  has Hausdorff
dimension $2$, and its
$2$--dimensional Hausdorff
measure is equal to $0$. 
\end{proposition}

\begin{proof} By
Corollary~\ref
{cor:IrrationalPoints}, it
suffices to show that $\mathbb
R-\mathbb Q$, endowed with the
metric $d$, has Hausdorff
dimension 2 and 2--dimensional
Hausdorff measure 0. By symmetry,
it suffices to show this for
$\left]0,\infty\right[ -\mathbb
Q$ endowed with the metric $d$.

We will make a further reduction.
The map $x \mapsto \frac x{x+1}$
induces an isometry between the
metric spaces
$\left( \left]0,\infty\right[
-\mathbb Q, d\right)$ and $\left(
\left[0,1\right]-\mathbb Q,
\delta \right)$, where 
\begin{equation*} \delta\bigl( x,y
\bigr) =\max
\bigl\{ {\textstyle \frac 1{q};
x\leq \frac{p}{q} \leq y }\bigr
\}.
\end{equation*} 
It therefore suffices to show
that $\left(
\left[0,1\right]-\mathbb Q,
\delta \right)$ has Hausdorff
dimension $2$ and 
$2$--dimensional Hausdorff
measure $0$. 

 We begin with a few elementary
observations on
$\delta$--balls in
$\left[0,1\right]-\mathbb Q$.
Recall that a \emph{Farey
interval} is an interval of the
form $\bigl[ \frac pq, \frac{p'}{q'}
\bigr]$ with $p'q-pq'=1$.

\begin{lemma} Let $B\subset
[0,1]-\mathbb Q$ be an open 
$\delta$--ball of radius $r$. Then
there is a Farey interval
$I=\bigl[ \frac pq, \frac{p'}{q'}
\bigr]$ such that $B = I-\mathbb
Q$. In addition, the diameter of
$(B,\delta)$ is equal to
$\Delta(B)=\Delta(I) =
\frac1{q+q'}$, and its Lebesgue
measure is $l(I)=\frac1{qq'}$.
\end{lemma}

\begin{proof} Let $N_r$ be the
finite set of all rational
numbers $\frac pq \in [0,1]$ with
$q\leq
\frac 1r$. It immediately follows
from the definition of the metric
$\delta$ that the ball $B$ is
equal to
$I-\mathbb Q$ for the closure
$I$ of some component of
$[0,1]-N_r$. It is well-known
that such an $I$ must be a Farey
interval; see for instance
\cite[Section~3]{HarWri}. The formula
for the diameter
$\Delta(B)$ is an immediate
consequence of the fact that
$\min \bigl\{ {q''}; \frac pq <
\frac{p''}{q''} < \frac {p'}{q'}
\bigr\}=q+q'$ for every  Farey
interval $I=\bigl[ \frac pq,
\frac{p'}{q'}
\bigr]$, which is elementary. The
formula for the length
$l(I)$ of the interval $I$ is a
straightforward computation.
\end{proof}

The reader should beware of
unexpected properties of
$\delta$--balls.  If $I=\bigl[
\frac pq, \frac{p'}{q'}
\bigr]$ is a Farey interval, then
$B=I-\mathbb Q$ is an open
$\delta$--ball whose center is any
point of
$B$ and whose radius is any
number $r$ with
$\frac1{q+q'}< r \leq \min \bigl\{
\frac1q,
\frac1{q'} \bigr\}$. 

\begin{lemma}
\label{lem:Hausdorffdim<2} The
Hausdorff dimension of
$\left( [0,1]-\mathbb Q,
\delta\right)$ is at least $2$.
\end{lemma}

\begin{proof} We will show that,
for every $s$ with $1<s<2$, the
$s$--dimensional Hausdorff
measure of
$\left[0,1\right]-\mathbb Q$ is
strictly positive.

 For this, we will use the
classical  fact
 that, for every
$\varepsilon>0$ and for almost
every $x\in
\left[0,1\right]$, the set
$\bigl\{ \frac pq \in \mathbb Q;
\bigl\vert x-\frac pq \bigr\vert
<\frac 1{q^{2+ \varepsilon}}
\bigr\}$ is finite; see for
instance
\cite{HarWri, Lang66}.  As a
consequence, there exists a subset
$A_\varepsilon\subset
\left[0,1
\right] - \mathbb Q$ with non-zero
Lebesgue measure and a number
$N_\varepsilon>0$ such that
$\bigl\vert x-\frac pq
\bigr\vert \geq\frac 1{q^{2+
\varepsilon}}$ for every $x \in
A_\varepsilon$ and every $\frac pq
\in
\mathbb Q$ with $q\geq
N_\varepsilon$.

We apply this to
$\varepsilon =\frac{2-s}{s-1}>0$.
If
$x\in A_\varepsilon$ is contained
in a Farey interval $I=\bigl[
\frac pq, \frac{p'}{q'} \bigr]$
where both $q$ and $q'$ are
greater then $ N_\varepsilon$,
suppose
$q'>q$ without loss of
generality. Then,
\begin{equation*}
\textstyle
\frac1{qq'} = l(I) 
\geq \left \vert x-\frac pq
\right\vert 
\geq \frac1{q^{2+\varepsilon}}
\end{equation*} so that $q\geq
\left(q'
\right)^{\frac1{1+\varepsilon}} =
\left( q' \right)^{s-1}$. It
follows that
\begin{equation*}
\textstyle
\Delta(I)^s = \frac1{\left( q+q'
\right)^{s}}
\geq \frac1{2^s \left( q'
\right)^{s}} > 
\frac1{4qq'} = \frac14 l(I)
\end{equation*} in this case.

Let $r>0$ be small. Cover
$\left[0,1\right]-
\mathbb Q$ by a family of
$\delta$--balls $B_i$, $i \in
\mathcal I$, with respective
radii $r_i \leq r$. By the
ultrametric property of $\delta$,
we can assume that the $B_i$ are
pairwise disjoint. (If two
$\delta$--balls meet, one is
contained in the other one). We
saw that each ball
$B_i$ is equal to $I_i-\mathbb Q$
for some Farey interval $I_i =
\bigl[ \frac{p_i}{q_i},
\frac{p'_i}{q'_i} \bigr]$. By the
ultrametric property, $r_i
\geq
\Delta(B_i) = \frac1{q_i +
q_i'}$. In particular, at least
one of $q_i$, $q_i'$ is greater
than
$\frac1{2r_i}\geq \frac1{2r}$, and
therefore is greater than
$N_\varepsilon$ if we choose
$r$ small enough. 

Decompose the index set $\mathcal
I$ as the disjoint union of
$\mathcal I_1$, $\mathcal I_2$
and $\mathcal I_3$ where:
\begin{enumerate}
\item $\mathcal I_1$ consists of
those $i$ such that $I_i$ contains
a point of
$A_\varepsilon$ and such that
both $q_i$ and $q_i'$ are greater
than $N_\varepsilon$.
\item $\mathcal I_2$ consists of
those $i$ such that $I_i$ meets
$A_\varepsilon$ and such that
only one of $q_i$ and $q_i'$ is
greater than $N_\varepsilon$.
\item $\mathcal I_3$ consists of
those $i$ such that $I_i$ does
not meet
$A_\varepsilon$.
\end{enumerate}

If $i\in \mathcal I_1$, namely if
$I_i$ contains a point of
$A_\varepsilon$ and if {both}
$q_i$ and $q_i'$ are greater than
$N_\varepsilon$, we saw that
$r_i^s \geq \Delta(B_i)^s
=\Delta(I_i)^s
\geq \frac14 l(I_i)$. 

If $i\in \mathcal I_2$, namely if
$I_i$ contains a point of
$A_\varepsilon$ and if only one
of $q_i$ and $q_i'$ is greater
than $N_\varepsilon$, we observed
that $\max \left\{ q_i, q_i'
\right\} \geq \frac 1{2r_i} \geq
\frac1{2r}$. Therefore, $l(I_i) =
\frac1{q_iq_i'} < 2r_i<2r$. Also,
note that there are at most
$N_\varepsilon^2$ such elements
of $\mathcal I_2$. 

We conclude that 
\begin{equation*} l(A_\varepsilon)
\leq \sum_{i\in \mathcal I_1}
l(I_i) +
\sum_{i\in \mathcal I_2} l(I_i)
\leq  4\sum_{i\in \mathcal I_1}
r_i^s + 2rN_\varepsilon^2
\leq  4\sum_{i\in \mathcal I}
r_i^s + 2rN_\varepsilon^2
\end{equation*}

Taking the infimum over all
coverings of $[0,1]-\mathbb Q$ by
$\delta$--balls of radius at most
$r$, and letting $r$ tend to 0, we
conclude that the $s$--dimensional
Hausdorff measure $\mathcal
H^s\left( [0,1]-\mathbb Q,
\delta\right) $ of $\left(
[0,1]-\mathbb Q, \delta\right)$ is
bounded from below by $\frac14
l(A_\varepsilon)>0$. 

This proves that $\left(
[0,1]-\mathbb Q, \delta\right)$
has non-zero
$s$--dimensional Hausdorff
measure for every $1<s<2$. It
follows that the Hausdorff
dimension of $\left(
[0,1]-\mathbb Q, \delta \right)$
is at least 2.
\end{proof}

\begin{lemma} The
$2$--dimensional Hausdorff
measure of $\left( [0,1]- \mathbb
Q, \delta\right)$ is equal to
$0$. In particular, its Hausdorff
dimension is at most $2$.
\end{lemma}
\begin{proof} We can use the
well-known connection between
Farey intervals and continued
fractions. If $x\in[0,1]-\mathbb
Q$ has continued fraction
expansion
$$x=\cfrac{1}{a_1(x)+ \cfrac{1}{a_2(x)
+ \cfrac{1}{
\cdots
+\cfrac1{a_n(x)
+\cfrac1{\cdots
}}}}}
=\left[ a_1(x), a_2(x),
\dots, a_n(x),
\dots\right],$$
 then successive
finite continued fractions
$\frac{p_n}{q_n}=
\left[ a_1(x), a_2(x),
\dots, a_n(x)\right]$ form Farey
intervals $\bigl[
\frac{p_n}{q_n},
\frac{p_{n-1}}{q_{n-1}} \bigr]$
(with $\frac{p_n}{q_n}>
\frac{p_{n-1}}{q_{n-1}}$ for $n$
odd) containing $x$. In addition,
$a_n (x)\leq
\frac{q_n}{q_{n-1}} < a_n(x) +1$
because
${q_{n}}= {a_n(x)q_{n-1}+q_{n-2}}$
and $0<q_{n-2}<q_{n-1}$.

The set $A$ of those 
$x=\left[ a_1(x), a_2(x),
\dots, a_n(x),
\dots\right]$ for which the
sequence $\left( a_n(x)\right)
_{n\in
\mathbb N}$ is unbounded has full
Lebesgue measure in $[0,1]$; see
for instance
\cite{HarWri, Lang66}. By
our comparison between 
$\frac{q_n}{q_{n-1}}$ and
$a_n(x)$, for every small $r>0$
and for every large
$M>1$, we can cover $A$ by a
family of Farey intervals
$I_i=\bigl[
\frac {p_i'}{ q_i'}, \frac
{p_i''}{q_i''} \bigr]$, $i\in
\mathcal I$, such that $
\frac {q_i'}{ q_i''} >M$ and
$q_i'+q_i'' >\frac1r$. (We do not
necessarily assume that 
$\frac {p_i'}{ q_i'} <\frac
{p_i''}{q_i''}$). For such an
interval, 
\begin{equation*}
\frac{\Delta(I_i)^2}{l(I_i)}
=\frac{q_i'  q_i''} {\left( q_i'
+q_i''\right)^2} = f
\left( \frac{q_i'}{q_i''} \right)
< f(M)
\end{equation*} where
$f(x)=
\frac{x}{\left(x+1\right)^2}$ and
where the inequality comes from
the fact that $f(x)$ is
decreasing for $x\geq 1$. In
particular, the intersection
$B_i$ of $I_i$ with
$[0,1]-\mathbb Q$ is a
$\delta$--ball of radius $r_i<r$
with $r_i^2 =\Delta(I_i)^2 <f(M)
l(I_i)$. 

Since the complement $B$ of $A$
has Lebesgue measure 0, we can
cover $B$ by a family of Farey
intervals $I_j$, $j\in \mathcal
J$, whose total length is
arbitrarily small, say $\sum_{j\in
\mathcal J} l(I_j) <\frac1M$. In
this case, the intersection
$B_j$ of $I_j$ with
$[0,1]-\mathbb Q$ is a $d$--ball
of radius $r_j<r$ with $r_j^2
=\Delta(I_j)^2\leq l(I_j)
\left(\max f(x) \right)= \frac14
l(I_j)$. 

In this way, we have covered
$[0,1]-\mathbb Q$ by a family of
$\delta$--balls $B_i$, $i\in
\mathcal I
\cup \mathcal J$, with respective
radii $r_i<r$ such that 
\begin{equation*}
\sum_{i\in \mathcal I \cup
\mathcal J} r_i^2 \leq f(M)
\sum_{i\in \mathcal I} l(I_i) +
{\textstyle\frac14} \sum_{j\in
\mathcal J}l(I_j)
\leq f(M) +
\frac1{4M}.
\end{equation*}

The estimate is valid for every
$r>0$ and every $M\geq 1$. Since
$f(M)$ tends to 0 as $M$ tends to
$\infty$, it follows that the
2--dimensional Hausdorff measure
of $\left( [0,1]- \mathbb
Q, \delta\right)$ is equal to
0. 
\end{proof}

This concludes the proof of
Proposition~\ref{prop:dim2}. 
\end{proof}

\Addresses
\recd
 
\end{document}